\documentclass[10pt,leqno]{article}
\usepackage[plainpages=false,pdfpagelabels]{hyperref}
\usepackage[final,notref,notcite]{showkeys}
\usepackage[a4paper,hscale=0.7,vscale=0.82,hmarginratio=1:1,includehead,headheight=14pt]{geometry}
\usepackage{fancyhdr}
\fancypagestyle{plain}{\fancyhead{}\fancyfoot{}}

\usepackage{amsthm}
\theoremstyle{definition}
\newtheorem{definition}{Definition}[section]
\newtheorem{remark}[definition]{Remark}
\newtheorem{example}[definition]{Example}
\newtheorem{problem}[definition]{Problem}
\theoremstyle{plain}
\newtheorem{lemma}[definition]{Lemma}
\newtheorem{proposition}[definition]{Proposition}
\newtheorem{theorem}[definition]{Theorem}
\newtheorem{corollary}[definition]{Corollary}
\newtheorem*{nntheorem}{Theorem}

\newenvironment{prf}{\begin{proof}[\textbf{\upshape Proof.}]}{\end{proof}}

\newenvironment{introduction}
	{\begin{center}{\textbf{\large Introduction}}\end{center}
	\markboth{Introduction}{Introduction}
}
	{}

\usepackage{amssymb,amsmath,amscd,dsfont}
\newcommand{\op}[1]{\ensuremath{\operatorname{#1}}}
\newcommand{\wt}[1]{\ensuremath{\widetilde{#1}}}
\newcommand{\wh}[1]{\ensuremath{\widehat{#1}}}
\newcommand{\ol}[1]{\ensuremath{\overline{#1}}}
\newcommand{\cl}[1]{\ensuremath{\overline{#1}}}
\newcommand{\fa}{\ensuremath{\;\text{\textup{ for all }}\;}}
\newcommand{\Exp}{\ensuremath{\operatorname{Exp}}}
\newcommand{\cE}{\ensuremath{\mathcal{E}}}

\newcommand{\cP}{\ensuremath{\mathcal{P}}}
\newcommand{\cV}{\ensuremath{\mathcal{V}}}

\newcommand{\fk}{\ensuremath{\mathfrak{k}}}
\newcommand{\fg}{\ensuremath{\mathfrak{g}}}

\newcommand{\1}{\ensuremath{\mathds{1}}}

\newcommand{\R}{\ensuremath{\mathds{R}}}
\newcommand{\N}{\ensuremath{\mathds{N}}}
\newcommand{\Z}{\ensuremath{\mathds{Z}}}

\newcommand{\K}{\ensuremath{\mathds{K}}}

\newcommand{\bS}{\ensuremath{\mathds{S}}}
\newcommand{\id}{\ensuremath{\operatorname{id}}}
\newcommand{\pr}{\ensuremath{\operatorname{pr}}}
\newcommand{\ev}{\ensuremath{\operatorname{ev}}}
\newcommand{\Ad}{\ensuremath{\operatorname{Ad}}}

\newcommand{\supp}{\ensuremath{\operatorname{supp}}}

\newcommand{\Gau}{\ensuremath{\operatorname{Gau}}}
\newcommand{\gau}{\ensuremath{\operatorname{\mathfrak{gau}}}}
\newcommand{\Aut}{\ensuremath{\operatorname{Aut}}}

\newcommand{\Diff}{\ensuremath{\operatorname{Diff}}}
\newcommand{\Conn}{\ensuremath{\operatorname{Conn}}}
\newcommand{\Gauc}{\ensuremath{\operatorname{Gau}_{c}}}
\newcommand{\Autc}{\ensuremath{\operatorname{Aut}_{c}}}
\newcommand{\Homeo}{\ensuremath{\operatorname{Homeo}}}
\newcommand{\Hom}{\ensuremath{\operatorname{Hom}}}

\newcommand{\glue}{\ensuremath{\operatorname{glue}}}
\newcommand{\res}{\ensuremath{\operatorname{res}}}

\newcommand{\In}{i=1,\ldots,n}

\newcommand{\Jm}{j=1,\ldots,m}

\newcommand{\pfb}[1][M]{\ensuremath{(K,\pi\from P \to #1)}}

\newcommand{\se}{\ensuremath{\nobreak\subseteq\nobreak}}
\newcommand{\from}{\ensuremath{\nobreak:\nobreak}}
\renewcommand{\to}{\ensuremath{\nobreak\rightarrow\nobreak}}

\newcommand{\diffeo}{g}
\newcommand{\form}{\omega}
\begin{document}
\title{\textbf{Lie Group Structures on Symmetry Groups of Principal Bundles}}
\author{Christoph Wockel\\
        Mathematisches Institut\\
        Georg-August-Universit\"at G\"ottingen\\
\small \texttt{christoph@wockel.eu}}
\maketitle
\thispagestyle{empty}
\begin{abstract}
In this paper we describe how one can obtain Lie group structures on
the group of (vertical) bundle automorphisms for a locally convex
principal bundle \mbox{$\cP$} over the compact manifold \mbox{$M$}. This is done by
first considering Lie group structures on the group of vertical bundle
automorphisms \mbox{$\Gau(\cP)$}. Then the full automorphism group
\mbox{$\Aut(\cP)$} is considered as an extension of the open subgroup
\mbox{$\Diff(M)_{\cP}$} of diffeomorphisms of \mbox{$M$} preserving the equivalence
class of \mbox{$\cP$} under pull-backs, by the gauge group \mbox{$\Gau(\cP)$}. We
derive explicit conditions for the extensions of these Lie group
structures, show the smoothness of some natural actions and relate our
results to affine Kac--Moody algebras and groups.\\[\baselineskip]
\\[\baselineskip]
\end{abstract}

\bigskip

\begin{introduction}

This paper enlarges the scope of infinite-dimensional Lie groups from
mapping groups to the group of bundle automorphisms and of vertical
bundle automorphisms (shortly called gauge group) of a principal
bundle \mbox{$\cP$} with compact base \mbox{$M$}. We do this by
introducing natural and easy accessible locally convex
Lie group structures on these groups, denoted by
\mbox{$\Aut(\cP)$} and \mbox{$\Gau(\cP)$}. They are interesting examples of
infinite-dimensional Lie groups, since they arise naturally as
symmetry groups of gauge field theories. In particular, it is shown
that \mbox{$\Aut(\cP)$} may also in a Lie theoretic context be
interpreted as an extension of (an open subgroup of) \mbox{$\Diff(M)$}
by \mbox{$\Gau(\cP)$}, such that it describes most naturally gauge
field symmetries and space-time symmetries in a unified way. Moreover,
for special types of bundles, \mbox{$\Gau(\cP)$} may be interpreted as
a Kac--Moody group and \mbox{$\Aut(\cP)$} as the automorphism group of
the corresponding Kac--Moody algebra (cf.\ Example
\ref{exmp:automorphismGroupOfTwistedLoopAlgebra}). This gives
geometric interpretations of these groups and leads to topological
information on them.

Lie group topologies on \mbox{$\Gau(\cP)$} and \mbox{$\Aut(\cP)$} have been
considered in the literature, e.g. in
\cite{omoriMaedaYoshiokaKobayashi83OnRegularFrechetLieGroups},
\cite{cirelleMania85theGroupOfGaugeTransformationsAsASchwartz-LieGroup},
\cite{abbati89},
\cite{michorAbbatiCirelliMania89automorphismGroup},
\cite{michor91GaugeTheoryForFiberBundles}, \cite{krieglmichor97},
mostly by methods of the Convenient Setting and for finite-dimensional
structure groups. These approaches have the common disadvantage that
they identify \mbox{$\Gau(\cP)$} and \mbox{$\Aut(\cP)$} with subsets of larger
manifolds of mappings with less structure, making it hard to access
the Lie group structure in concrete computations. Unfortunately, the
proof given in \cite[Theorem 42.21]{krieglmichor97} has a serious gap
since the chart used in the construction needs to be defined on an
\mbox{$\Ad$}-invariant zero neighbourhoods in the Lie algebra of the
structure group, which does not exist the generality claimed
there. The key advantage of the approach taken in this paper is that
we provide explicit charts for the Lie group structures without any
compactness assumptions on the structure group, making it possible to
check smoothness conditions directly, as illustrated in case of some
natural actions in Proposition
\ref{prop:actionOfAutomorphismGroupOnBundleIsSmooth} and Proposition
\ref{lem:automorphismGroupActingOnConnections}. In addition, the power
of the approach to the Lie group structure on \mbox{$\Gau(\cP)$} from Section
\ref{sect:topologisationOfTheGaugeGroup} becomes clear in Section
\ref{sect:theFullAutomorphismGroupAsInfiniteDimensionalLieGroup},
where smoothness conditions for \mbox{$\Gau(\cP)$}-valued cocycles are
explicitly checked. However, it should be emphasised that our interest
in these groups is not a gauge-theoretic one but comes from
infinite-dimensional Lie theory. This will also become clear from the
applications we provide.

We now describe our results in some detail. In the first section we
introduce the Lie group structure on \mbox{$\Gau(\cP)$}, quite similar to the
case of mapping groups \mbox{$C^{\infty}(M,K)$} as in \cite{pressleysegal86}
or \cite{gloeckner02b}, but there is a little subtlety. In order to
make this approach work we have to impose a technical condition on
\mbox{$\cP$}, ensuring the compatibility of charts of the structure group \mbox{$K$}
and transition functions of the bundle \mbox{$\cP$}. This condition is called
``property SUB'' and will ensure the constructions to work throughout
the whole paper. Our fist result is then the following.

\begin{nntheorem}[Lie group structure on \mbox{$\boldsymbol{\Gau(\cP)}$}] Let
\mbox{$\cP$} be a smooth principal \mbox{$K$}-bundle over the compact manifold \mbox{$M$}
(possibly with corners).  If \mbox{$\cP$} has the property SUB, then the
gauge group \mbox{\mbox{$\Gau(\cP)\cong C^{\infty}(P,K)^{K}$}} carries a Lie
group structure, modelled on \mbox{$C^{\infty}(P,\fk)^{K}$}.  If, moreover,
\mbox{$K$} is locally exponential , then \mbox{$\Gau(\cP)$} is so.
\end{nntheorem}

In the remainder of the section we elaborate on the property SUB and,
moreover, show that a broad variety of bundles have these
property. Since all prominent classes of infinite-dimensional Lie
groups occur as structure groups of bundles having the property SUB it
is a natural question whether this is the case for all locally convex
Lie groups but this question remains unanswered.
\vskip\baselineskip

In the second section we enlarge the scope of Lie group structures to
the full automorphism group \mbox{$\Aut(\cP)$}. This is done by considering
\mbox{$\Aut(\cP)$} as an extension of the open subgroup \mbox{$\Diff(M)_{\cP}$} of
diffeomorphisms of \mbox{$M$} which preserve the equivalence class of \mbox{$\cP$}
under pull-backs, by \mbox{$\Gau(\cP)$}. Using the recently established
theory of non-abelian Lie group extensions for locally convex Lie
groups from \cite{neeb06nonAbelianExtensions} we obtain our second
result.

\begin{nntheorem}[\mbox{$\boldsymbol{\Aut(\cP)}$} as an extension of
\mbox{$\boldsymbol{\Diff(M)_{\cP}}$} by \mbox{$\boldsymbol{\Gau(\cP)}$}]
Let \mbox{$\cP$} be a smooth
principal \mbox{$K$}-bundle over the closed compact manifold \mbox{$M$}. If \mbox{$\cP$}
has the property SUB, then \mbox{$\Aut(\cP)$} carries a Lie group structure
such that we have an extension of smooth Lie groups
\begin{align*}
\Gau(\cP)\hookrightarrow \Aut(\cP)\xrightarrow{Q}
\Diff(M)_{\cP},
\end{align*}
where \mbox{$Q:\Aut(\cP)\to \Diff(M)$} is the natural homomorphism and
\mbox{$\Diff(M)_{\cP}$} is the open subgroup of \mbox{$\Diff(M)$} preserving the
equivalence class of \mbox{$\cP$} under pull-backs.
\end{nntheorem}

An interesting thing about this theorem is that it relates two major
classes of infinite-dimensional Lie group in a non-trivial way, namely
mapping groups, having many ideals from evaluation homomorphisms, and
groups of diffeomorphisms, which are perfect (cf.\
\cite{hallerTeichmann03SmoothPerfectnessThroughDecomposingDiffeomorphisms}
and references therein). Moreover, two points are important about the
proof of this theorem. The first is that it avoids connections and
differential equations and thus also works beyond regularity, taking
the property SUB as sole technical requirement. Furthermore, this
procedure provides explicit charts for the Lie group structure, making
it possible to prove smoothness of natural actions as already
mentioned.

In the rest of the section we illustrate the interpretation of
\mbox{$\Aut(\cP)$} as automorphism group of certain types of Kac--Moody
Algebras. As an immediate consequence of the preceding theorem we
obtain a long exact homotopy sequence, explaining some of the
topological properties of \mbox{$\Aut(\cP)$}.
\end{introduction}

\section{The Lie group topology on the gauge group}
\label{sect:topologisationOfTheGaugeGroup}
We shall mostly identify the gauge group with the space of
\mbox{$K$}-equivariant continuous mappings \mbox{$C^{\infty}(P,K)^{K}$}, where \mbox{$K$}
acts on itself by conjugation from the right.
This identification allows us to topologise the gauge group very
similar to mapping groups \mbox{$C^{\infty}(M,K)$} for compact \mbox{$M$}. Since the
compactness of \mbox{$M$} is the crucial point in the topologisation of
mapping groups, we can not take this approach directly, because our
structure groups \mbox{$K$} shall not be compact, even infinite-dimensional.

\begin{definition}\label{def:bundleAutomorphismsAndGaugeTransformations}
If \mbox{$K$} is a topological group and \mbox{$\cP=\pfb$} is a continuous
principal \mbox{$K$}-bundle, then we denote by \[
\Autc(\cP):=\{f\in\Homeo
(P):\rho_{k}\op{\circ}f=f\op{\circ}\rho_{k}\fa k\in K\}
\]
the
group of continuous \emph{bundle automorphisms}  and by \[
\Gauc(\cP):=\{f\in \Autc(\cP):\pi \op{\circ}f=\pi \}
\]
the group of continuous \emph{vertical} bundle automorphisms or
\emph{continuous gauge group}. If, in addition, \mbox{$K$} is a Lie group,
\mbox{$M$} is a manifold with corners and \mbox{$\cP$} is a smooth principal bundle,
then we denote by \[ \Aut(\cP):=\{f\in\Diff
(P):\rho_{k}\op{\circ}f=f\op{\circ}\rho_{k}\fa k\in K\}
\]
the the
group of \emph{smooth bundle automorphisms}  (or shortly
\emph{bundle automorphisms}). Then
each \mbox{$F\in \Aut(\cP)$} induces an element \mbox{$F_{M}\in \Diff
(M)$} , given
by \mbox{$F_{M}(p\cdot K):=F(p)\cdot K$} if we identify \mbox{$M$} with \mbox{$P/K$}. This
yields a homomorphism \mbox{$Q\from \Aut(\cP)\to \Diff(M)$}, \mbox{$F\mapsto F_{M}$}
and we denote by \mbox{$\Gau(\cP)$} the kernel of \mbox{$Q$} and by \mbox{$\Diff
(M)_{\cP}$} the image of \mbox{$Q$}. Thus \[
\Gau(\cP)=\{f\in \Aut(\cP): \pi \op{\circ}f=\pi\},
\]
which we call the group of (smooth) vertical bundle
automorphisms or shortly the \emph{gauge group}  of \mbox{$\cP$}.
\end{definition}

\begin{remark}\label{rem:gaugeGroupIsIsomorphicToEquivariantMappings}
If \mbox{$\cP$} is a smooth principal \mbox{$K$}-bundle and if we denote by \[
C^{\infty}(P,K)^{K}:=\{\gamma \in C^{\infty}(P,K):
\gamma (p\cdot k)= k^{-1}\cdot \gamma (p)\cdot k\fa p\in P,k\in K\}
\] the group of \mbox{$K$}-equivariant smooth maps from \mbox{$P$} to \mbox{$K$}, then the map
\[
C^{\infty}(P,K)^{K}\ni f\mapsto \big(p\mapsto p\cdot
f(p)\big)\in\Gau(\cP)
\]
is an isomorphism   of groups and we will
mostly identify \mbox{$\Gau(\cP)$} with \mbox{$C^{\infty}(P,K)^{K}$} via this map.
\end{remark}

The algebraic counterpart of the gauge group is the gauge
algebra. This will serve as the modelling space for the gauge group
later on.

\begin{definition}
If \mbox{$\cP$} is a smooth principal \mbox{$K$}-bundle, then the space \[
\gau(\cP):=C^{\infty}(P,\fk)^{K}:=\{\eta \in C^{\infty}(P,\fk )^{K}:
\eta (p\cdot k)=\Ad(k^{-1}).\eta (p)\fa p\in P,k\in K \}
\]
is called the \emph{gauge algebra}  of \mbox{$\cP$}. We endow it with the subspace
topology from \mbox{$C^{\infty}(P,\fk)$} and with the pointwise Lie
bracket.
\end{definition}

\begin{proposition}\label{prop:isomorphismOfTheGaugeAlgebra} Let
\mbox{$\cP=\pfb$} be a smooth principal \mbox{$K$}-bundle over the
finite-di\-men\-sion\-al manifold with corners \mbox{$M$}. If \mbox{$\cl{\cV}
:=(\cl{V}_{i},\sigma_{i})_{i\in I}$} is a smooth closed trivialising
system of \mbox{$\cP$} with transition functions \mbox{$k_{ij}\from \cl{V}_{i}\cap
\cl{V}_{j}\to K$}, then we denote
\[
\fg_{\cl{\cV}}(\cP):=\left\{(\eta_{i})_{i\in I}\in\prod_{i\in I}
C^{\infty}(\cl{V}_{i},\fk ): \eta_{i}(m)=\Ad
(k_{ij}(m)).\eta_{j}(m)\fa m\in \cl{V}_{i}\cap \cl{V}_{j} \right\}.
\]
If \mbox{$\cV$} denotes the smooth open trivialising system underlying
\mbox{$\cl{\cV}$}, then we set \[
\fg_{\cV}(\cP):=\left\{(\eta_{i})_{i\in I}\in\prod_{i\in I}
C^{\infty}(V_{i},\fk ):
\eta_{i}(m)=\Ad(k_{ij}(m)).\eta_{j}(m)\fa m\in V_{i}\cap V_{j}\right\},
\]
and we have isomorphisms of topological vector spaces
\[
\gau(\cP)= C^{\infty}(P,\fk)^{K}
\cong S(\Ad(\cP))\cong \fg_{\cV}(\cP)\cong\fg_{\ol{\cV}}(\cP).
\]
Furthermore, each of these spaces is a locally convex Lie algebra in a
natural way and the isomorphisms are isomorphisms of topological Lie
algebras.
\end{proposition}

\begin{prf}
The last two isomorphisms are provided by Proposition
\ref{prop:gluingForVectorBundles_CLOSED_Version} and Corollary
\ref{cor:gluingForVectorBundles_OPEN_Version}, so we show
\mbox{$C^{\infty}(P,\fk)^{K}\cong \fg_{\cl{\cV}}(\cP)$}. 

For each \mbox{$\eta\in C^{\infty}(P,\fk)^{K}$} the element \mbox{$(\eta_{i})_{i\in
I}$} with \mbox{$\eta_{i}=\eta \circ \sigma_{i}$}
defines an element of \mbox{$\fg_{\cl{\cV}}(\cP)$} and the map
\[
\psi \from C^{\infty}(P,\fk)^{K}\to \fg_{\cl{\cV}}(\cP),\quad
\eta\mapsto(\eta_{i})_{i\in I}
\]
is continuous. In fact, \mbox{$\sigma_{i}(m)=\sigma_{j}(m)\cdot k_{ji}(m)$} for
\mbox{$m\in \cl{V}_{i}\cap \cl{V}_{j}$} implies
\[
\eta_{i}(m)=\eta (\sigma_{i}(m))=\eta (\sigma_{j}(m)\cdot k_{ji}(m))
=\Ad (k_{ji}(m))^{-1}.\eta(\sigma_{j}(m))=\Ad(k_{ij}(m)).\eta_{j}(m)
\]
and thus \mbox{$(\eta_{i})_{i\in I}\in\fg_{\cl{\cV}}(\cP)$}. Recall that if
\mbox{$X$} is a topological space, then a map \mbox{$f\from X\to
C^{\infty}(\cl{V}_{i},\fk)$} is continuous if and only if \mbox{$x\mapsto
d^{n}f(x)$} is continuous for each \mbox{$n\in \N_{0}$} (Remark
\ref{rem:alternativeDescriptionOfTopologyOnSpaceOfSmoothMappings}). This
implies that \mbox{$\psi$} is continuous, because
\mbox{$d^{n}\eta_{i}=d^{n}\eta\op{\circ}T^{n}\sigma_{i}$} and pull-backs
along continuous maps are continuous.

On the other hand, if \mbox{$k_{i}\from \pi^{-1}(\cl{V}_{i})\to K$} is
given by \mbox{$p=\sigma_{i}(\pi (p))\cdot k_{i}(p)$} and if \mbox{$(\eta_{i})_{i\in
I}\in \fg_{\cl{\cV}}(\cP)$}, then the map
\[
\eta:P\to \fk,\quad p\mapsto
\Ad\left(k(p)\right)^{-1}.\eta_{i}\left(\pi (p)\right) \;\text{ if }\;
\pi (p)\in \cl{V}_{i}
\]
is well-defined, smooth and \mbox{$K$}-equivariant. Furthermore,
\mbox{$(\eta_{i})_{i\in I}\mapsto \eta$} is an inverse of \mbox{$\psi$} and it thus
remains to check that it is continuous, i.e., that
\[
\fg_{\cl{\cV}}(\cP)\ni(\eta_{i})_{i\in I}\mapsto d^{n}\eta \in C(T^{n}P,\fk)
\]
is continuous for all \mbox{$n\in \N_{0}$}. If \mbox{$C\se T^{n}P$} is compact, then
\mbox{$(T^{n}\pi) (C)\se T^{n}M$} is compact and hence it is covered by
finitely many \mbox{$T^{n}V_{i_{1}},\dots ,T^{n}V_{i_{m}}$} and thus
\mbox{$\left(T^{n}\left(\pi^{-1}(\ol{V_{i}})\right)\right)_{i=i_{1},\dots
,i_{m}}$} is a finite closed cover of \mbox{$C\se T^{n}P$}. Hence it suffices
to show that the map
\[
\fg_{\cl{\cV}}(\cP)\ni(\eta_{i})_{i\in I}\mapsto T^{n}(\left.\eta
\right|_{\pi^{-1}(\cl{V}_{i})}) \in C(T^{n}\pi^{-1}(\cl{V}_{i}),\fk )
\]
is continuous for \mbox{$n\in \N_{0}$} and \mbox{$i\in I$} and we may thus
w.l.o.g. assume that \mbox{$\cP$} is trivial. In the trivial case we have
\mbox{$\eta=\Ad(k^{-1}).(\eta\circ \pi)$} if \mbox{$p\mapsto (\pi (p),k(p))$}
defines a global trivialisation.  We shall make the case \mbox{$n=1$}
explicit. The other cases can be treated similarly and since the
formulae get quite long we skip them here.

Given any open zero neighbourhood in \mbox{$C(TP,\fk )$}, which we may assume
to be \mbox{$\lfloor C,V\rfloor$} with \mbox{$C\se TP$} compact and \mbox{$0\in V\se \fk$}
open, we have to construct an open zero neighbourhood \mbox{$O$} in
\mbox{$C^{\infty}(M,\fk)$} such that \mbox{$\varphi(O)\se \lfloor C,V\rfloor$}.  For
\mbox{$\eta '\in C^{\infty}(M,\fk)$} and \mbox{$X_{p}\in C$} we get with Lemma
\ref{lem:productrule}
\[
d(\varphi(\eta ' ))(X_{p})= \Ad(k^{-1}(p)).d\eta'(T\pi(X_{p}))-
[\delta^{l}(k)(X_{p}),\Ad(k^{-1}(p)).\eta' (\pi(p)) ].
\]
Since \mbox{$\delta^{l}(C)\se \fk$} is compact, there exists an open zero
neighbourhood \mbox{$V'\se \fk$} such that
\[
\Ad(k^{-1}(p)).V'+[\delta^{l}(k)(X_{p}),\Ad(k^{-1}(p)).V']\se
V
\]
for each \mbox{$X_{p}\in C$}. Since \mbox{$T\pi :TP\to TM$} is continuous, \mbox{$T\pi
(C)$} is compact and we may set \mbox{$O=\lfloor T\pi (C),V'\rfloor$}.

That \mbox{$\fg_{\cV}(\cP)$} and \mbox{$\fg_{\cl{\cV}}(\cP)$} are locally convex Lie
algebras follows because they are closed subalgebras of \mbox{$\prod_{i\in
I}C^{\infty}(V_{i},\fk)$} and \mbox{$\prod_{i\in I}C^{\infty}(\cl{V}_{i},\fk
)$}. Since the isomorphisms
\[
C^{\infty}(P,\fk)^{K}
\cong S(\Ad(\cP))\cong \fg_{\cV}(\cP)\cong\fg_{\ol{\cV}}(\cP).
\]
are all isomorphisms of abstract Lie algebras an isomorphisms of
locally convex vector spaces, it follows that they are isomorphisms
of topological Lie algebras.
\end{prf}

\begin{definition}\label{def:localGaugeGroup}
If \mbox{$\cP$} is a smooth \mbox{$K$}-principal bundle with compact base and
\mbox{$\cl{\cV}=(\cl{V}_{i},\sigma_{i})_{\In}$} is a smooth closed
trivialising system with corresponding transition functions
\mbox{$k_{ij}\from\cl{V}_{i}\cap \cl{V}_{j}\to K$}, then we denote \[
G_{\ol{\cV }}(\cP):=\left\{(\gamma_{i})_{\In}\in
\prod_{i=1}^{n}C^{\infty}(\ol{V_{i}},K):
\gamma_{i}(m)=k_{ij}(m)\gamma_{j}k_{ji}(m) \fa m\in \cl{V}_{i}\cap
\cl{V}_{j}\right\}
\]
and turn it into a group with respect to pointwise group operations.
\end{definition}

\begin{remark}\label{rem:isoToTheGaugeGroupInLocalCoordinates} In the
situation of Definition \ref{def:localGaugeGroup}, the map
\begin{align}
\label{eqn:localGaugeGroupIsAbstractlyIsomorphisToEquivariantMappings}
\psi \from G_{\cl{\cV}}(\cP)\to C^{\infty}(P,K)^{K},\quad \psi
((\gamma_{i})_{\In})(p)= k^{-1}_{\sigma_{i}}(p)\cdot
\gamma_{i}(\pi(p))\cdot k_{\sigma_{i}}(p)\;\text{ if }\; \pi (p)\in
\cl{V}_{i}
\end{align}
is an isomorphism of abstract groups, where the map on the right hand
side is well-defined because \mbox{$k_{\sigma_{i}}(p)=k_{ij}(\pi (p))\cdot
k_{\sigma_{j}}(p)$} and thus
\begin{multline*}
k_{\sigma_{i}}^{-1}(p)\cdot \gamma_{i}(\pi (p))\cdot k_{\sigma_{i}}(p)=
k_{\sigma_{j}}(p)^{-1}\cdot 
\underbrace{k_{ji}(\pi (p))\cdot \gamma_{i}(\pi (p))\cdot k_{ij}(\pi (p))}%
_{\gamma_{j}(\pi (p))}\cdot k_{\sigma_{j}(p)}\\
=k_{\sigma_{j}}(p)^{-1}\cdot \gamma_{j}(\pi (p))\cdot k_{\sigma_{j}}(p).
\end{multline*}
In particular, this implies that \mbox{$\psi ((\gamma_{i})_{\In})$} is
smooth. Since for \mbox{$m\in \ol{V}_{i}$} the evaluation map
\mbox{$\ev_{m}:C^{\infty}(\ol{V}_{i},K)\to K$} is continuous,
\mbox{$G_{\cl{\cV}}(\cP)$} is a closed subgroup of
\mbox{$\prod_{i=1}^{n}C^{\infty}(\ol{V}_{i},K)$}.
\end{remark}

Since an infinite-dimensional Lie group may posses closed subgroups
which are no Lie groups (cf.\ \cite[Exercise III.8.2]{bourbakiLie}),
the preceding remark does not automatically yield a Lie group
structure on \mbox{$G_{\cl{\cV}}(\cP)$}.  However, in many situations, it
will turn out that \mbox{$G_{\cl{\cV}}(\cP)$} has a natural Lie group
structure.

The following definition encodes the necessary requirement ensuring a
Lie group structure on \mbox{$G_{\cl{\cV}}(\cP)$} that is induced by the
natural Lie group structure on
\mbox{$\prod_{i=1}^{n}C^{\infty}(\cl{V}_{i},K)$}. Since quite different
properties of \mbox{$\cP$} will ensure this requirement it seems to be worth
extracting it as a condition on \mbox{$\cP$}. The name for this requirement
will be justified in Corollary \ref{cor:propertySUB}. 

\begin{definition}\label{def:propertySUB}
If \mbox{$\cP$} is a smooth principal \mbox{$K$}-bundle with compact base and 
\mbox{$\cl{\cV}=(\cl{V}_{i},\sigma_{i})_{\In}$} is a smooth closed
trivialising system, then we say that \mbox{$\cP$} has the \emph{property
SUB} with respect to \mbox{$\cl{\cV}$} if there exists a
convex centred chart \mbox{$\varphi\from W\to W'$} of \mbox{$K$} such that
\[
\varphi_{*}\from G_{\cl{\cV}}(\cP)\cap \prod_{i=1}^{n}C^{\infty}(\cl{V}_{i},W)
\to \fg_{\cl{\cV}}(\cP)\cap \prod_{i=1}^{n}C^{\infty }(\cl{V}_{i},W'),
\quad (\gamma_{i})_{\In}\mapsto (\varphi \op{\circ}\gamma_{i})_{\In}
\]
is  bijective. We say that \mbox{$\cP$} has the property SUB if \mbox{$\cP$}
has this property with respect to some trivialising system.
\end{definition}

It should be emphasised that in all relevant cases, known to the
author, the bundles have the property SUB, and it is still unclear,
whether there are bundles, which do not have this property (cf.\ Lemma
\ref{lem:propertySUB} and Remark \ref{rem:propertySUB}). This property
now ensures the existence of a natural Lie group structure on
\mbox{$G_{\cl{\cV}}(\cP)$}.

\begin{proposition}\label{prop:gaugeGroupInLocalCoordinatesIsLieGroup}
\textup{\textbf{a)}} Let \mbox{$\cP$} be a smooth principal \mbox{$K$}-bundle with
compact base \mbox{$M$}, which has the property SUB with respect to the
smooth closed trivialising system \mbox{$\cl{\cV}$}. Then \mbox{$\varphi_{*}$}
induces a smooth manifold structure on \mbox{$G_{\cl{\cV}}(\cP)\cap
\prod_{i=1}^{n}C(\cl{V}_{i},W)$}. Furthermore, the conditions \mbox{$i)-iii)$}
of Proposition \ref{prop:localDescriptionsOfLieGroups} are satisfied
such that \mbox{$G_{\cl{\cV}}(\cP)$} can be turned into a Lie group modelled
on \mbox{$\fg_{\cl{\cV}}(\cP)$}.

\textup{\textbf{b)}} In the setting of \textup{\textbf{a)}}, the map \mbox{$\psi
\from G_{\cl{\cV}}(\cP)\to C^{\infty}(P,K)^{K}$} is an isomorphism of
topological groups if \mbox{$C^{\infty}(P,K)^{K}$} is endowed with the
subspace topology from \mbox{$C^{\infty}(P,K)$}.

\textup{\textbf{c)}} In the setting of \textup{\textbf{a)}}, we have
\mbox{$\op{L}(G_{\cl{\cV }}(\cP))\cong\fg_{\cl{\cV }}(\cP)$}.
\end{proposition}

\begin{prf}
\textbf{a)} Set \mbox{$U:=G_{\cl{\cV}}(\cP)\cap
\prod_{i=1}^{n}C(\cl{V}_{i},K)$}.  Since \mbox{$\varphi_{*}$} is bijective by
assumption and \mbox{$\varphi_{*}(U)$} is open in \mbox{$\fg_{\cl{\cV}}(\cP)$}, it induces
a smooth manifold structure on \mbox{$U$}.

Let \mbox{$W_{0}\se W$} be an open unit neighbourhood with \mbox{$W_{0}\cdot
W_{0}\se W$} and \mbox{$W_{0}^{-1}=W_{0}$}. Then \mbox{$U_{0}:=G_{\cl{\cV }}(\cP)\cap
\prod_{i=1}^{n}C^{\infty}(\ol{V_{i}},W_{0})$} is an open unit
neighbourhood in \mbox{$U$} with \mbox{$U_{0}\cdot U_{0}\se U$} and
\mbox{$U_{0}=U_{0}^{-1}$}.  Since each \mbox{$C^{\infty}(\ol{V_{i}},K)$} is a
topological group, there exist for each \mbox{$(\gamma_{i})_{\In}$} open unit
neighbourhoods \mbox{$U_{i}\se C^{\infty}(\ol{V_{i}},K)$} with
\mbox{$\gamma_{i}\cdot U_{i}\cdot \gamma_{i}^{-1}\se
C^{\infty}(\ol{V_{i}},W)$}. Since \mbox{$C^{\infty}(\ol{V_{i}},W_{0})$} is
open in \mbox{$C^{\infty}(\ol{V_{i}},K)$}, so is \mbox{$U'_{i}:=U_{i}\cap
C^{\infty}(\ol{V_{i}},W_{0})$}.  Hence
\[
(\gamma_{i})_{\In} \cdot \left(G_{\cl{\cV}}(\cP)\cap(U'_{1}\times \dots \times
U'_{n}) \right)\cdot (\gamma_{i}^{-1})_{\In}\se U
\]
and conditions \mbox{$i)-iii)$} of Proposition
\ref{prop:localDescriptionsOfLieGroups} are satisfied, where the
required smoothness properties are consequences of
the smoothness of push forwards of mappings between function spaces
(cf.\ \cite[Proposition 28 and Corollary 29]{smoothExt} and
\cite[3.2]{gloeckner02a}).

\textbf{b)} We show that the map \mbox{$\psi \from G_{\cl{\cV}}(\cP)\to
C^{\infty}(P,K)^{K}$} from
\eqref{eqn:localGaugeGroupIsAbstractlyIsomorphisToEquivariantMappings}
is a homeomorphism. Let \mbox{$\left.\cP\right|_{\cl{V}_{i}}=:\cP_{i}$} be
the restricted bundle. Since \mbox{$T^{n}\cl{V}_{i}$} is closed in \mbox{$T^{n}M$},
we have that \mbox{$C^{\infty}(P,K)^{K}$} is homeomorphic to
\[
\wt{G}_{\cl{\cV}}(\cP):=\{(\wt{\gamma}_{i})_{\In}\in 
\prod_{i=1}^{n}C^{\infty}(P_{i},K)^{K}:\wt{\gamma}_{i}(p)=\wt{\gamma}_{j}(p)
\fa p\in \pi^{-1}(\cl{V}_{i}\cap \cl{V}_{j})\}
\]
as in Proposition \ref{prop:gluingForVectorBundles_CLOSED_Version}.
With respect to this identification, \mbox{$\psi$} is given by
\[
(\gamma_{i})_{\In}\mapsto (k_{\sigma_{i}}^{-1}\cdot (\gamma_{i}\op{\circ}\pi) 
\cdot k_{\sigma_{i}})_{\In}
\]
and it thus suffices to show the assertion for trivial bundles.
So let \mbox{$\sigma \from M\to P$} be a global section.
The map \mbox{$C^{\infty}(M,K)\ni f\mapsto f\circ \pi \in C^{\infty}(P,K)$}
is continuous since
\[
C^{\infty}(M,K)\ni f\mapsto T^{k}(f\circ \pi)=T^{k}f\circ T^{k}\pi =
(T^{k}\pi)_{*}(T^{k}f)\in C(T^{k}P,T^{k}K)
\]
is continuous as a composition of a pullback an the map \mbox{$f\mapsto
T^{k}f$}, which defines the topology on \mbox{$C^{\infty}(M,K)$}. Since
conjugation in \mbox{$C^{\infty}(P,K)$} is continuous, it follows that
\mbox{$\varphi$} is continuous.  Since the map \mbox{$f\mapsto f\circ \sigma$} is
also continuous (with the same argument), the assertion follows.

\textbf{c)} This follows immediately from
\mbox{$\op{L}(C^{\infty}(\ol{V_{i}},K))\cong C^{\infty}(\ol{V_{i}},\fk)$}
(cf.\ \cite[Section 3.2]{gloeckner02a}).
\end{prf}

The next corollary is a mere observation. Since it justifies the 
name ``property SUB'', it is made explicit here.

\begin{corollary}\label{cor:propertySUB}
If \mbox{$\cP$} is a smooth principal \mbox{$K$}-bundle with compact base \mbox{$M$},
having the property SUB with respect to the smooth closed trivialising
system \mbox{$\cl{\cV}$}, then \mbox{$G_{\cl{\cV}}(\cP)$} is a closed subgroup of
\mbox{$\prod_{i=1}^{n}C^{\infty}(\cl{V}_{i},K)$}, which is a Lie group
modelled on \mbox{$\fg_{\cl{\cV}}(\cP)$}.
\end{corollary}

That different choices of charts lead to isomorphic Lie group
structures follows directly from Proposition
\ref{prop:localDescriptionsOfLieGroups}. We show next that in fact,
different choices of trivialising systems also lead to isomorphic Lie
group structures on \mbox{$\Gau(\cP)$}.

\begin{proposition}
\label{prop:isomorhicLieGroupsStructuresOnGau(P)ForDifferentTrivialisations}
Let \mbox{$\cP$} be a smooth principal \mbox{$K$}-bundle with compact
base. If we have two trivialising systems
\mbox{$\cl{\cV}=(\cl{V}_{i},\sigma_{i})_{\In}$} and
\mbox{$\cl{\mathcal{U}}=(\cl{U}_{j},\tau_{j})_{\Jm}$} and \mbox{$\cP$}
has the property SUB with respect to \mbox{$\cl{\cV}$} and
\mbox{$\cl{\mathcal{U}}$}, then \mbox{$G_{\ol{\cV}}(\cP)$} is
isomorphic to \mbox{$G_{\ol{\mathcal{U}}}(\cP)$} as a Lie group.
\end{proposition}

\begin{prf}
First, we note that if the covers underlying \mbox{$\cl{\cV}$} and
\mbox{$\cl{\mathcal{U}}$} are the same, but the sections differ by smooth
functions \mbox{$k_{i}\in C^{\infty}(\cl{V}_{i},K)$}, i.e.,
\mbox{$\sigma_{i}=\tau_{i}\cdot k_{i}$}, then this induces an automorphism
of Lie groups
\[
G_{\cl{\cV}}(\cP)\to G_{\cl{\cV}}(\cP),\quad 
(\gamma_{i})_{\In}\mapsto (k_{i}^{-1}\cdot\gamma_{i}\cdot k_{i})_{\In},
\]
because conjugation with \mbox{$k_{i}^{-1}$} is an automorphism of
\mbox{$C^{\infty}(\cl{V}_{i},K)$}.

Since each two open covers have a common refinement it suffices to
show the assertion if one cover is a refinement of the other.  So let
\mbox{$V_{1},\dots ,V_{n}$} be a refinement of \mbox{$U_{1},\dots ,U_{m}$} and let
\mbox{$\{1,\dots ,n\} \ni i\mapsto j(i)\in \{1,\dots ,m\}$} be a function
with \mbox{$V_{i}\se U_{j(i)}$}. Since different choices of sections lead to
automorphisms we may assume that
\mbox{$\sigma_{i}=\left.\sigma_{j(i)}\right|_{\cl{V}_{i}}$}, implying in
particular \mbox{$k_{ii'}(m)=k_{j(i)j(i')}(m)$}. Then the restriction map
from Lemma \ref{lem:restrictionMapForCurrentGroupIsSmooth} yields a
smooth homomorphism
\[
\psi :G_{\cl{\mathcal{U}}}(\cP)\to G_{\cl{\cV}}(\cP),\quad
(\gamma_{j})_{j\in J}\mapsto
(\left.\gamma_{j(i)}\right|_{\cl{V}_{i}})_{i\in I}.
\]

For \mbox{$\psi^{-1}$} we construct each component \mbox{$\psi_{j}^{-1}\from
G_{\cl{\cV}}(\cP)\to C^{\infty}(\cl{U}_{j},K)$} separately. The condition that
\mbox{$(\psi_{j}^{-1})_{j\in J}$} is inverse to \mbox{$\psi$} is then
\begin{align}
\label{eqn:isomorhicLieGroupsStructuresOnGau(P)ForDifferentTrivialisations1} 
\left.\psi^{-1}_{j}((\gamma_i)_{i\in I})\right|_{\cl{V}_{i}}=\gamma_{i}
\fa i\;\text{ with }\; j=j(i).
\end{align}

Set \mbox{$I_{j}:=\{i\in I: \cl{V}_{i}\se \cl{U}_{j}\}$} and note that
\mbox{$j(i)=j$} implies \mbox{$i\in I_{j}$}.  Since a change of the sections
\mbox{$\sigma_{i}$} induces an automorphism on \mbox{$G_{\cl{\cV}}(\cP)$} we may
assume that \mbox{$\sigma_{i}=\left.\sigma_{j(i)}\right|_{\cl{V}_{i}}$} for
each \mbox{$i\in I_{j}$}.  Let \mbox{$x \in \cl{U}_{j}\backslash \cup_{i\in
I_{j}}V_{i}$}. Then \mbox{$x\in V_{i_{x}}$} for some \mbox{$i_{x}\in I$} and thus
there exists an open neighbourhood \mbox{$U_{x}$} of \mbox{$x$} such that
\mbox{$\cl{U}_{x}$} is a manifold with corners, contained in \mbox{$\cl{U}_{j}\cap
\cl{V}_{i_{x}}$}.  Now finitely many
\mbox{$U_{x_{1}},\dots ,U_{x_{l}}$} cover \mbox{$\cl{U}_{j}\backslash \cup_{i\in
I_{j}}V_{i}$} and we set
\[
\psi^{-1}_{j}((\gamma_{i})_{i\in I})=\glue \left((\gamma_{i})_{i\in
I_{j}}, \left(\left.(k_{ji_{x_{k}}}\cdot\gamma_{i_{x_{k}}}\cdot
k_{i_{x_{k}}j})\right|_{U_{x_{k}}} \right)_{k=1,\dots ,l} \right).
\]
Then this defines a smooth map by Proposition
\ref{prop:gluingLemmaForCurrentGroup} and
\eqref{eqn:isomorhicLieGroupsStructuresOnGau(P)ForDifferentTrivialisations1}
is satisfied because \mbox{$j(i)=i$} implies \mbox{$i\in I_{j}$}
\end{prf}

We now come to the main result of this section.

\begin{theorem}[Lie group structure on \mbox{$\boldsymbol{\Gau(\cP)}$}]
\label{thm:gaugeGroupIsLieGroup} Let \mbox{$\cP$} be a smooth principal
\mbox{$K$}-bundle over the compact manifold \mbox{$M$} (possibly with corners).  If
\mbox{$\cP$} has the property SUB, then the gauge group \mbox{\mbox{$\Gau(\cP)\cong
C^{\infty}(P,K)^{K}$}} carries a Lie group structure, modelled on
\mbox{$C^{\infty}(P,\fk)^{K}$}.  If, moreover, \mbox{$K$} is locally exponential ,
then \mbox{$\Gau(\cP)$} is so.
\end{theorem}

\begin{prf}
We endow \mbox{$\Gau(\cP)$} with the Lie group structure induced from the
isomorphisms of groups \mbox{$\Gau(\cP)\cong C^{\infty}(P,K)^{K}\cong
G_{\cl{\cV}}(\cP)$} for some smooth closed trivialising system
\mbox{$\cl{\cV}$}. To show that \mbox{$\Gau(\cP)$} is locally exponential if \mbox{$K$} is
so we first show that if \mbox{$M$} is a compact manifold with corners and
\mbox{$K$} has an exponential function, then
\[
(\exp_{K})_{*}:C^{\infty}(M,\fk )\to C^{\infty}(M,K),\quad \eta \mapsto
\exp_{K}\op{\circ }\eta 
\]
is an exponential function for \mbox{$C^{\infty}(M,K)$}. For \mbox{$x\in \fk$} let
\mbox{$\gamma_{x}\in C^{\infty}(\R,K)$} be the solution of the initial value
problem \mbox{$\gamma (0)=e$}, \mbox{$\gamma' (t)=\gamma (t).x$}. Take \mbox{$\eta \in
C^{\infty}(M,\fk)$}.  Then
\[
\Gamma_{\eta}\from \R\to C^{\infty}(M,K),\quad (t,m)\mapsto \gamma_{\eta (m)}(t)
=\exp_{K}(t\cdot \eta (m))
\]
is a homomorphism of abstract groups. Furthermore, \mbox{$\Gamma_{\eta}$} is
smooth, because it is smooth on a zero neighbourhood of \mbox{$\R$}, for the
push-forward of the local inverse of \mbox{$\exp_{K}$} provide charts on a
unit neighbourhood in \mbox{$C^{\infty}(M,K)$}. Then
\[
\delta^{l}(\Gamma_{\eta})(t)=\Gamma_{\eta}(t)^{-1}\cdot \Gamma
'(t)=\Gamma_{\eta}(t)^{-1}\cdot \Gamma_{\eta}(t)\cdot \eta=\eta,
\]
thought of as an equation in the Lie group \mbox{$T\big(C^{\infty}(M,K)\big)\cong
C^{\infty}(M,\fk)\rtimes C^{\infty}(M,K)$}, shows that \mbox{$\eta \mapsto
\Gamma_{\eta}(1)=\exp_{K}\circ \gamma$} is an exponential function for
\mbox{$C^{\infty}(M,K)$}. The proof of the preceding lemma yields immediately
that
\[
\fg_{\cl{\cV}}(\cP)\cap
\prod_{i=1}^{n}C^{\infty}(\ol{V_{i}},W')\to
G_{\cl{\cV}}(\cP),\quad
(\eta_{i})_{\In}\mapsto
(\exp_{K}\circ \eta )_{\In}
\]
is a diffeomorphism and thus \mbox{$\Gau(\cP)$} is locally exponential.
\end{prf}

It remains to elaborate on the arcane property SUB. First we shall see
that this property behaves well with respect to refinements of
trivialising systems.

\begin{lemma}\label{lem:propertySUBIsCompatibleWithRefinements} Let
\mbox{$\cP$} be a smooth principal \mbox{$K$}-bundle with compact base and
\mbox{$\cl{\cV}=(\cl{V}_{i},\sigma_{i})_{\In}$} be a smooth closed
trivialising system of \mbox{$\cP$}. If
\mbox{$\cl{\mathcal{U}}=(\cl{U}_{j},\tau_{j})_{\Jm }$} is a refinement of
\mbox{$\cl{\cV}$}, then \mbox{$\cP$} has the property SUB with respect to \mbox{$\cl{\cV}$}
if and only if \mbox{$\cP$} has the property SUB with respect to
\mbox{$\cl{\mathcal{U}}$}.
\end{lemma}

\begin{prf}
Let \mbox{$\{1,\dots ,m\}\ni j\mapsto i(j)\in\{1,\dots ,n\}$} be a map such
that \mbox{$U_{j}\se V_{i(j)}$} and
\mbox{$\tau_{j}=\left.\sigma_{i(j)}\right|_{\cl{U}_{j}}$}. Then we have
bijective mappings
\begin{alignat*}{2}
\psi_{G}&\from G_{\cl{\cV}}(\cP)\to G_{\cl{\mathcal{U}}}(\cP),\quad 
&&(\gamma_{i})_{\In}\mapsto (\left.\gamma_{i(j)}\right|_{\cl{j}})_{\Jm}\\
\psi_{\fg}&\from \fg_{\cl{\cV}}(\cP)\to \fg_{\cl{\mathcal{U}}}(\cP),\quad 
&&(\eta_{i})_{\In}\mapsto (\left.\eta_{i(j)}\right|_{\cl{j}})_{\Jm}
\end{alignat*}
(cf.\ Proposition
\ref{prop:isomorhicLieGroupsStructuresOnGau(P)ForDifferentTrivialisations}).
Now let \mbox{$\varphi \from W\to W'$} be an arbitrary convex centred chart
of \mbox{$K$} and set
\begin{gather*}
Q:=G_{\cl{\cV}}(\cP)\cap \prod_{i=1}^{n}C(\cl{V}_{i},W)\quad\quad
\wt{Q}:=G_{\cl{\mathcal{U}}}(\cP)\cap \prod_{i=1}^{n}C(\cl{U}_{i},W)\\
Q':=\fg_{\cl{\cV}}(\cP)\cap \prod_{i=1}^{n}C(\cl{V}_{i},W')\quad\quad
\wt{Q}':=\fg_{\cl{\mathcal{U}}}(\cP)\cap \prod_{i=1}^{n}C(\cl{U}_{i},W')
\end{gather*}
Then we have \mbox{$\psi_{G}(Q)=\wt{Q}$} and \mbox{$\psi_{\fg}(Q')=\wt{Q}'$} and the
assertion follows from the commutative diagram
\[
\begin{CD}
Q@>\varphi_{*}>>Q'\\
@V\psi_{G}VV @V\psi_{\fg}VV\\
\wt{Q}@>\varphi_{*}>>\wt{Q}'.
\end{CD}
\]
\end{prf}

The following lemma will be needed later in the proof that bundles
with structure group a direct limit Lie group have the property SUB.

\begin{lemma}\label{lem:equivariantExtensionOfCharts}
Let \mbox{$C$} be a compact Lie group, \mbox{$M,N$} be smooth finite-dimensional
manifolds with \mbox{$N\se M$} and assume that the inclusion \mbox{$i\from
M\hookrightarrow N$} is a smooth immersion. Furthermore, let \mbox{$C$} act on
\mbox{$N$} from the right such that \mbox{$M\cdot C=M$} and that \mbox{$x\in M$} is
fix-point. Then \mbox{$C$} acts on \mbox{$T_{x}N$}, leaving \mbox{$T_{x}M$} invariant.

If there exists a \mbox{$C$}-invariant smoothly and \mbox{$C$}-equivariantly
contractible relatively compact open neighbourhood \mbox{$U$} of \mbox{$x$} and a
\mbox{$C$}-equivariant chart \mbox{$\varphi \from U\to \wt{U}\se T_{x}M$} with
\mbox{$\varphi (x)=0$}, then there exists a \mbox{$C$}-equivariant chart \mbox{$\psi \from
V\to \wt{V}\se T_{x}N$} such that \mbox{$V$} is a \mbox{$C$}-invariant
smoothly and \mbox{$C$}-equivariantly contractible relatively compact open
neighbourhood of \mbox{$x$} in \mbox{$M$}, satisfying \mbox{$V\cap M=U$}, \mbox{$\psi (V)\cap
T_{x}M=\varphi (U)$} and \mbox{$\left.\psi\right|_{U}=\varphi$}.
\end{lemma}

\begin{prf}(cf. \cite[Lemma 2.1]{gloeckner05} for the non-equivariant
case) Fix a \mbox{$C$}-invariant metric on \mbox{$TN$}, inducing a \mbox{$C$}-invariant
metric on \mbox{$N$} (c.f.\ \cite[Section
VI.2]{bredon72transformationGroups}). Then \mbox{$C$} acts on \mbox{$N$} by
isometries. As in \cite[Lemma 2.1]{gloeckner05} we find a
\mbox{$\sigma$}-compact, relatively compact open submanifold \mbox{$W'$} of \mbox{$N$} such
that \mbox{$W'\cap \ol{U}=U$}, whence \mbox{$U$} is a closed submanifold of \mbox{$W'$}.  We
now set \mbox{$W:=\cup_{c\in C}W'\cdot c$}.  Since \mbox{$C$} acts by isometries,
this is an open \mbox{$C$}-invariant subset of \mbox{$N$} and we deduce that we
still have
\[
W\cap \ol{U}=\left(\cup_{c\in C}W'\cdot c\right)\cap \left(\cup_{c\in
C}\ol{U}\cdot c \right)=\cup_{c\in C}\left((W'\cap \ol{U})\cdot c\right)=U.
\]
By shrinking \mbox{$W$} if necessary we thus get an open \mbox{$C$}-invariant
relatively compact submanifold of \mbox{$M$} with \mbox{$U$} as closed submanifold.

By \cite[Theorem VI.2.2]{bredon72transformationGroups}, \mbox{$U$} has an
open \mbox{$C$}-invariant tubular neighbourhood in \mbox{$W$}, i.e., there exists a
\mbox{$C$}-vector bundle \mbox{$\xi \from E\to U$} and a \mbox{$C$}-equivariant
diffeomorphism \mbox{$\Phi \from E\to W$} onto some \mbox{$C$}-invariant open
neighbourhood \mbox{$\Phi (E)$} of \mbox{$W$} such that the restriction of \mbox{$\Phi$} to
the zero section \mbox{$U$} is the inclusion of \mbox{$U$} in \mbox{$W$}. The proof of
\cite[Theorem 2.2]{bredon72transformationGroups} shows that \mbox{$E$} can be
taken to be the normal bundle of \mbox{$U$} with the canonical \mbox{$C$}-action,
which is canonically isomorphic to the \mbox{$C$}-invariant subbundle
\mbox{$TU^{\perp}$}. We will thus identify \mbox{$E$} with \mbox{$TU^{\perp}$} from now on.

That \mbox{$U$} is smoothly and \mbox{$C$}-equivariantly contractible
means that there exists a homotopy \mbox{$F\from [0,1]\times U\to U$}
such that each \mbox{$F(t,\cdot )\from U\to U$} is smooth and
\mbox{$C$}-equivariant and that \mbox{$F(1,\cdot)$} is the map which
is constantly \mbox{$x$} and \mbox{$F(0,\cdot )=\id_{U}$}. Pulling
back \mbox{$E$} along the smooth and equivariant map
\mbox{$F(1,\cdot)$} gives the \mbox{$C$}-vector bundle
\mbox{$\pr_{1}\from U\times T_{x}U^{\perp}\to U$}, where the action of
\mbox{$C$} on \mbox{$U$} is the one given by assumption and the action
of \mbox{$C$} on \mbox{$T_{x}U^{\perp}$} is the one induced from the
canonical action of \mbox{$C$} on \mbox{$T_{x}M$}. By \cite[Corollary
2.5]{wasserman69EquivariantDifferentialTopology}, \mbox{$F(1,\cdot
)^{*}(TU^{\perp})$} and \mbox{$F(0,\cdot)^{*}(TU^{\perp})$} are
equivalent \mbox{$C$}-vector bundles and thus there exists a smooth
\mbox{$C$}-equivariant bundle equivalence \mbox{$\Psi \from
TU^{\perp}=F(1,\cdot )^{*}(TU^{\perp})\to U\times
T_{x}U^{\perp}=F(0,\cdot)^{*}(TU^{\perp})$}.

We now define
\[
\psi\from V:=\Phi (TU^{\perp})\to T_{x}N,\quad y\mapsto 
\varphi \left(\Psi_{1}(\Phi^{-1}(y))\right)+\Psi_{2}\left(\Phi^{-1}(y)\right),
\]
where \mbox{$\Psi_{1}$} and \mbox{$\Psi_{2}$} are the components of \mbox{$\Psi$}. Since
\mbox{$\Phi$}, \mbox{$\Psi$} and \mbox{$\varphi$} are \mbox{$C$}-equivariant so is \mbox{$\psi$} and it
is a diffeomorphism onto the open subset \mbox{$U\times T_{x}U^{\perp}$}
because \mbox{$T_{x}N=T_{x}M\oplus T_{x}U^{\perp}$}. This yields a
\mbox{$C$}-equivariant chart. Moreover, \mbox{$V$} is relatively compact as a subset
of \mbox{$W$}. Furthermore, if we denote by \mbox{$\psi_{1}$} the \mbox{$T_{x}M$}-component
and by \mbox{$\psi_{2}$} the \mbox{$T_{x}U^{\perp}$}-component of \mbox{$\psi$}, then
\[
[0,1]\times V\to V,\quad (t,y)\mapsto (\psi_{1}(y),t\cdot \psi_{2}(y))
\]
defines a smooth \mbox{$C$}-equivariant homotopy from \mbox{$\id_{V}$} to the
projection \mbox{$V\to U$}. Composing this homotopy with the smooth
\mbox{$C$}-equivariant homotopy from \mbox{$\id_{U}$} to the map which is constantly
\mbox{$x$} yields the asserted homotopy. Since \mbox{$V\cap M$} is the zero section
in \mbox{$\left.TU^{\perp}\right|_{U}$}, we have \mbox{$V\cap M=U$}. Furthermore, we
have \mbox{$\psi (V)\cap T_{x}M=\varphi (U)$}, because \mbox{$\Phi$} and \mbox{$\Psi$}
restrict to the identity on the zero section. This also implies
\mbox{$\left.\psi\right|_{U}=\varphi $}.  We thus have checked all
requirements from the assertion.
\end{prf}

Although it is presently unclear, which bundles have the property SUB
and which not, we shall now see that \mbox{$\cP$} has the property SUB in
many interesting cases, covering large classes of presently known
locally convex Lie groups.

\begin{lemma}\label{lem:propertySUB}
Let \mbox{$\cP$} be a smooth principal \mbox{$K$}-bundle over the compact manifold
with corners \mbox{$M$}.
\begin{itemize}
\item [\textup{\textbf{a)}}] If \mbox{$\cP$} is trivial, then there exists
a global smooth trivialising system and \mbox{$\cP$} has the
property SUB with respect to each such system.

\item [\textup{\textbf{b)}}] If \mbox{$K$} is abelian, then \mbox{$\cP$} has the property
SUB with respect to each smooth closed trivialising system.

\item [\textup{\textbf{c)}}] If \mbox{$K$} is a Banach--Lie group, then \mbox{$\cP$} has the
property SUB with respect to each smooth closed trivialising system.

\item [\textup{\textbf{d)}}] If \mbox{$K$} is locally exponential, then \mbox{$\cP$} has the
property SUB with respect to each smooth closed trivialising system.

\item [\textup{\textbf{e)}}] If \mbox{$K$} is a countable direct limit of
finite-dimensional Lie groups in the sense of \cite{gloeckner05}, then
there exists a smooth closed trivialising system such that the
corresponding transition functions take values in a compact subgroup
\mbox{$C$} of some \mbox{$K_{i}$} and \mbox{$\cP$} has the property SUB with respect to
each such system.
\end{itemize}
\end{lemma}

\begin{prf}
\textup{\textbf{a)}} If \mbox{$\cP$} is trivial, then there exists a global section 
	\mbox{$\sigma \from M\to P$} and thus \mbox{$\cl{\cV}=(M,\sigma)$} is a trivialising
	system of \mbox{$\cP$}. Then \mbox{$G_{\cl{\cV}}(\cP)=C^{\infty}(M,K)$} and
	\mbox{$\varphi_{*}$} is bijective for any convex centred chart 
	\mbox{$\varphi \from W\to W'$}.

\textup{\textbf{b)}} 
	If \mbox{$K$} is abelian, then the conjugation action of \mbox{$K$} on itself
	and the adjoint action of \mbox{$K$} on \mbox{$\fk$} are trivial. Then a direct
	verification shows that \mbox{$\varphi_{*}$} is bijective for any
	trivialising system \mbox{$\cl{\cV}$} and any convex centred chart
        \mbox{$\varphi \from W\to W'$}.

\textup{\textbf{c)}} If \mbox{$K$} is a Banach--Lie group, then it is in particular
	locally exponential (cf.\ Remark 
	\ref{rem:banachLieGroupsAreLocallyExponential})
	and it thus suffices to show \textup{\textbf{d)}}.

\textup{\textbf{d)}} Let \mbox{$K$} be locally exponential and
	\mbox{$\cl{\cV}=(\cl{V}_{i},\sigma_{i})_{\In}$} be a trivialising system.
	Furthermore, let \mbox{$W'\se \fk$} be an open zero neighbourhood such
	that \mbox{$\exp_{K}$} restricts to a diffeomorphism on \mbox{$W'$} and set
	\mbox{$W=\exp (W')$} and \mbox{$\varphi:=\exp^{-1} \from W\to W'$}.
	Then we have
\[
(\gamma_{i})_{\In}\in G_{\cl{\cV}}(\cP)\cap \prod_{i=1}^{n}C(\cl{V}_{i},W)
\;\;\Leftrightarrow\;\;
\varphi_{*}((\gamma_{i})_{\In})
\in \fg_{\cl{\cV}}(\cP)\cap \prod_{i=1}^{n}C(\cl{V}_{i},W'),
\]
because \mbox{$\exp_{K}(\Ad(k).x)=k\cdot
\exp_{K}(x)\cdot k^{-1}$} holds for all \mbox{$k\in K$} and \mbox{$x\in W'$}
(cf.\ Lemma \ref{lem:interchangeOfActionsOnGroupAndAlgebra}).
Furthermore, \mbox{$(\eta_{i})_{\In}\mapsto (\exp\op{\circ}\eta_{i})_{\In}$}
provides an inverse to \mbox{$\varphi_{*}$}.

\textup{\textbf{e)}} Let \mbox{$K$} be a direct limit of the countable direct
system \mbox{$\mathcal{S}=\left((K_{i})_{i\in I},(\lambda_{i,j})_{i\geq j}
\right)$} of finite-dimensional Lie groups \mbox{$K_{i}$} and Lie group
morphisms \mbox{$\lambda_{i,j}\from K_{j}\to K_{i}$} with
\mbox{$\lambda_{i,j}\op{\circ}\lambda_{j,\ell}=\lambda_{i,\ell}$} if \mbox{$i\geq
j\geq \ell$}. Then there exists an associated injective quotient system
\mbox{$\left((\ol{K}_{i})_{i\in I},(\ol{\lambda}_{i,j})_{i\geq j}\right)$}
with \mbox{$\ol{K}_{i}=K_{i}/N_{i}$}, where \mbox{$N_{i}=\ol{\bigcup_{j\geq i}\ker
\lambda_{j,i}}$} and \mbox{$\ol{\lambda}_{i,j}\from \ol{K}_{j}\to \ol{K}_{i}$}
is determined by
\mbox{$q_{i}\op{\circ}\lambda_{i,j}=\ol{\lambda}_{i,j}\op{\circ}q_{j}$} for
the quotient map \mbox{$q_{i}\from K_{i}\to \ol{K}_{i}$}. In particular, each
\mbox{$\ol{\lambda}_{i,j}$} is an injective immersion.

After passing to a cofinal subsequence of an equivalent direct system
(cf.\ \cite[\S{}1.6]{gloeckner05}), we may without loss of generality
assume that \mbox{$I=\N$}, that \mbox{$K_{1}\se K_{2}\se \dots$} and that the
immersions are the inclusion maps. Then a chart \mbox{$\varphi\from W\to W'$}
of \mbox{$K$} around \mbox{$e$} is the direct limit of a sequence of charts
\mbox{$(\varphi_{i}\from W_{i}\to W'_{i})_{i\in \N}$} such that
\mbox{$W=\bigcup_{i\in \N}W_{i}$}, \mbox{$W'=\bigcup_{i\in \N}W_{i}'$} and that
\mbox{$\left.\varphi_{i}\right|_{W_{j}}=\varphi_{j}$} if \mbox{$i\geq j$} (cf.\
\cite[Theorem 3.1]{gloeckner05}).

Now, let \mbox{$\cl{\cV}=(\ol{V}_{i},\ol{\sigma}_{i})_{\In}$} be a smooth
closed trivialising system of \mbox{$\cP$}. Then the corresponding transition
functions are defined on the compact subset \mbox{$\ol{V}_{i}\cap
\ol{V}_{j}$} and thus take values in a compact subset of \mbox{$K$}. Since
each compact subset of \mbox{$K$} is entirely contained in one of the \mbox{$K_{i}$}
(cf.\ \cite[Lemma 1.7]{gloeckner05}), the transition functions take values in
some \mbox{$K_{a}$}. Since each finite-dimensional principal \mbox{$K_{a}$}-bundle
can be reduced to a \mbox{$C$}-bundle, where \mbox{$C\se K_{a}$} is the maximal
compact Subgroup of \mbox{$K_{a}$} (cf.\ \cite[p.\ 59]{steenrod51}), we find
smooth mappings \mbox{$f_{i}\from \ol{V}_{i}\to K_{i}$} such that
\mbox{$\ol{\tau}_{i}:=\ol{\sigma}_{i}\cdot f_{i}$} is a smooth closed
trivialising system of \mbox{$\cP$} and the corresponding transition
functions take values in the compact Lie group \mbox{$C$}.

We now define a chart \mbox{$\varphi \from W\to W'$} satisfying the
requirements of Definition \ref{def:propertySUB}. For \mbox{$i<a$} let
\mbox{$W_{i}$} and \mbox{$W_{i}'$} be empty. For \mbox{$i=a$} denote
\mbox{$\fk_{a}:=\op{L}(K_{a})$} and let \mbox{$\exp_{a}\from \fk_{a}\to K_{a}$} be
the exponential function of \mbox{$K_{a}$}. Now \mbox{$C$} acts on \mbox{$K_{a}$} from the
right by conjugation and on \mbox{$\fk_{a}$} by the adjoint representation,
which is simply the induced action on \mbox{$T_{e}K_{a}$} for the fixed-point
\mbox{$e\in K_{a}$}. By Lemma \ref{lem:interchangeOfActionsOnGroupAndAlgebra}
we have \mbox{$\exp_{a}(\Ad (c).x)=c\cdot \exp_{a}(x)\cdot c^{-1}$} for
each \mbox{$c\in C$}. We choose a \mbox{$C$}-invariant metric on \mbox{$\fk_{a}$}. Then
there exists an \mbox{$\varepsilon >0$} such that \mbox{$\exp_{a}$} restricts to a
diffeomorphism on the open \mbox{$\varepsilon$}-ball \mbox{$W_{a}'$} around
\mbox{$0\in\fk_{a}$} in the chosen invariant metric. Then
\[
\varphi_{a}:=(\left.\exp_{a}\right|_{W'_{a}})^{-1}\from
W_{a}:=\exp_{a}(W_{a}')\to W'_{a},\quad \exp_{a}(x)\mapsto x
\]
defines an equivariant chart of \mbox{$K_{a}$} for the corresponding
\mbox{$C$}-actions on \mbox{$K_{a}$} and \mbox{$\fk_{a}$}, and, moreover, we may choose
\mbox{$W_{a}'$} so that \mbox{$\exp_{a}(W'_{a})$} is relatively compact in
\mbox{$K_{a}$}. In addition,
\[
[0,1]\times W_{a}\ni(t,k)\mapsto \exp_{a}\left(t\cdot \varphi (k)
\right)\in W_{a}
\]
defines a smooth \mbox{$C$}-equivariant contraction of \mbox{$W_{a}$}. By Lemma
\ref{lem:equivariantExtensionOfCharts} we may extend \mbox{$\varphi_{a}$} to
an equivariant chart \mbox{$\varphi_{a+1}\from W_{a+1}\to W_{a+1}'$} with
\mbox{$W_{a+1}\cap K_{a}=W_{a}$}
\mbox{$\left.\varphi_{a+1}\right|_{W_{a}}=\varphi_{a}$} such that \mbox{$W_{a+1}$}
is relatively compact in \mbox{$K_{a+1}$} and smoothly and \mbox{$G$}-equivariantly
contractible. Proceeding in this way we define \mbox{$G$}-equivariant charts
\mbox{$\varphi_{i}$} for \mbox{$i\geq a$}. 

This yields a direct limit chart \mbox{$\varphi :=\lim_{\to}\varphi_{i}\from
W\to W'$} of \mbox{$K$} for which we have \mbox{$W=\bigcup_{i\in I}W_{i}$} and \mbox{$W'=\bigcup_{i\in
I}W'_{i}$}. Since the action of \mbox{$C$} on \mbox{$T_{e}K_{i}=\fk_{i}$} is the
induced action in each step and the construction yields
\mbox{$\varphi_{i}(c^{-1}\cdot k\cdot c)=\Ad (c^{-1}).\varphi_{i}(k)$} we
conclude that we have
\[
\varphi^{-1}(\Ad (c^{-1}).x)=c^{-1}\cdot \varphi^{-1}(x)\cdot c\fa
x\in W'\text{ and }
c\in C
\]
(note that \mbox{$\exp$} is \emph{not} an inverse to \mbox{$\varphi$} any
more). Since the transition functions of the trivialising system
\mbox{$\cl{\cV}':=(\ol{\tau}_{i},\ol{V}_{i})$} take values in \mbox{$C$}, we may
proceed as in \textup{\textbf{d)}} to see that \mbox{$\varphi_{*}$} is
bijective.
\end{prf}

\begin{remark}\label{rem:propertySUB}
The preceding lemma shows that there are different kinds of properties
of \mbox{$\cP$} that can ensure the property SUB, i.e., topological in case
\textbf{a)}, algebraical in case \textbf{b)} and geometrical in case
\textbf{d)}. Case \textbf{e)} is even more remarkable, since it
provides examples of principal bundle with the property SUB, whose
structure groups are not locally exponential in general (c.f.\
\cite[Remark 4.7]{gloeckner05}).  It thus seems to be hard to find a
bundle which does \emph{not} have this property. However, a more
systematic answer to the question which bundles have this property is
not available at the moment.
\end{remark}

\begin{problem}
Is there a smooth principal \mbox{$K$}-bundle \mbox{$\cP$} over a compact base space \mbox{$M$}
which does not have the property SUB?
\end{problem} Lie group structures on the gauge group have already been considered
by other authors in similar settings.

\begin{remark}
If the structure group \mbox{$K$} is the group of diffeomorphisms \mbox{$\Diff(N)$}
of some closed compact manifold \mbox{$N$}, then it does not follow from
Lemma \ref{lem:propertySUB} that \mbox{$\cP$} has the property SUB, because
\mbox{$\Diff(N)$} fails to be locally exponential or abelian. However, in
this case, \mbox{$\Gau(\cP)$} is as a split submanifold of the Lie group
\mbox{$\Diff(P)$}, which provides a smooth structure on \mbox{$\Gau(\cP)$}
\cite[Theorem 14.4]{michor91GaugeTheoryForFiberBundles}.

Identifying \mbox{$\Gau(\cP)$} with the space of section in the
associated bundle \mbox{$\op{AD}(\cP)$} for the conjugation action
\mbox{$\op{AD}:K\times K\to K$}, \cite[Proposition
6.6]{omoriMaedaYoshiokaKobayashi83OnRegularFrechetLieGroups} also
provides a Lie group structure on \mbox{$\Gau(\cP)$}.

The advantage of Theorem \ref{thm:gaugeGroupIsLieGroup} is, that it
provides charts for \mbox{$\Gau(\cP)$}, which allows us to reduce questions
on gauge groups to similar question on mapping groups.  This
correspondence is crucial for all the following considerations.
\end{remark}

In the end of the section we provide as approximation result that
makes homotopy groups of \mbox{$\Gau(\cP)$} accessible in terms of continuous
data from the bundle \mbox{$\cP$} (cf.\ \cite{connHom}).

\begin{remark}\label{rem:topologyOnContinuousGaugeGroup}
Let \mbox{$\cP$} be one of the bundles that occur in Lemma
\ref{lem:propertySUB} with the corresponding closed trivialising
system \mbox{$\cl{\cV}$} and let \mbox{$\fg_{\cl{\cV}}(\cP)_{c}$} (resp.\
\mbox{$G_{\cl{\cV}}(\cP )_{c}$}) be the continuous counterparts of
\mbox{$\fg_{\cl{\cV}}(\cP)$} (resp.\ \mbox{$G_{\cl{\cV}}$}), which is isomorphic to
the space of \mbox{$K$}-equivariant continuous maps \mbox{$C(P,\fk)^{K}$} (resp.\
\mbox{$C(P,K)^{K}$}). We endow all spaces of continuous mappings with the
compact-open topology. Then the map
\[
\varphi_{*}\from G_{\cl{\cV}}(\cP)_{c}\cap \prod_{i=1}^{n}C(\cl{V}_{i},W)
\to \fg_{\cl{\cV}}(\cP)_{c}\cap \prod_{i=1}^{n}C(\cl{V}_{i},W'),
\quad (\gamma_{i})_{\In}\mapsto (\varphi \op{\circ}\gamma_{i})_{\In}
\]
is also bijective, inducing a smooth manifold structure on the
left-hand-side, because the right-hand-side is an open subset in a
locally convex space. Furthermore, it can be shown exactly as in the
smooth case, that this manifold structure endows
\mbox{$G_{\cl{\cV}}(\cP)_{c}$} with a Lie group structure by Proposition
\ref{prop:localDescriptionsOfLieGroups}. One could go on and call the
requirement that \mbox{$\varphi_{*}$} is bijective ``\textit{continuous
property SUB}'', but since the next proposition is the sole
application of it this seems to be exaggerated.
\end{remark}

\begin{lemma}\label{lem:approximationLemma2}
Let \mbox{$\cP$} be a smooth principal \mbox{$K$}-bundle over the compact base
\mbox{$M$}, having the property SUB with respect to the smooth closed
trivialising system \mbox{$\cl{\cV}=(\cl{V}_{i},\sigma_{i})_{\In}$} and let
\mbox{$\varphi \from W\to W'$} be the corresponding chart of
\mbox{$K$} (cf.\ Definition \ref{def:propertySUB}). If \mbox{$(\gamma_{i})_{\In}\in
G_{\cl{\cV}}(\cP)$} represents an element of \mbox{$C^{\infty}(P,K)^{K}$}
(cf.\ Remark \ref{rem:isoToTheGaugeGroupInLocalCoordinates}), which is
close to identity, in the sense that \mbox{$\gamma_{i} (\ol{V_{i}})\se W$},
then \mbox{$(\gamma_{i})_{\In}$} is homotopic to the constant map \mbox{$(x\mapsto
e)_{\In}$}.
\end{lemma}

\begin{prf}
Since the map
\[
\varphi_{*}:U:=G_{\cl{\cV}}(\cP)\cap \prod_{i=1}^{n}C^{\infty}(\ol{V_{i}},W)\to
\fg (\cP),\;\; (\gamma'_{i})_{\In} \mapsto
(\varphi \op{\circ}\gamma'_{i})_{\In},
\]
is a chart of \mbox{$G_{\cl{\cV}}(\cP)$} (cf.\ Proposition
\ref{prop:gaugeGroupInLocalCoordinatesIsLieGroup}) and
\mbox{$\varphi_{*}(U)\se \fg_{\cl{\cV}}(\cP)$} is convex, the map
\[
[0,1]\ni t \mapsto \varphi_{*}^{-1}\big(t\cdot
\varphi_{*}((\gamma_{i})_{\In})\big)\in G_{\cl{\cV}}(\cP)
\]
defines the desired homotopy.
\end{prf}

\begin{proposition}\label{thm:weakHomotopyEquivalence}
If \mbox{$\cP$} is one of the bundles that occur in Lemma
\ref{lem:propertySUB}, the natural inclusion \mbox{$\iota \from\Gau(\cP)
\hookrightarrow \Gauc(\cP)$} of smooth into continuous gauge
transformations is a weak homotopy equivalence, i.e., the induced
mappings \mbox{$\pi_{n}(\Gau(\cP))\to\pi_{n}\left(\Gauc(\cP)\right)$} are
isomorphisms of groups for \mbox{$n\in \N_{0}$}.
\end{proposition}

\begin{prf}
We identify \mbox{$\Gau(\cP)$} with \mbox{$C^{\infty}(P,K)^{K}$} and \mbox{$\Gauc(\cP)$}
with \mbox{$C(P,K)^{K}$}.  To see that \mbox{$\pi_{n}(\iota)$} is surjective,
consider the continuous principal \mbox{$K$}-bundle \mbox{$\pr^{*}(\cP)$} obtained
form \mbox{$\cP$} by pulling it back along the projection \mbox{$\pr:\bS^{n}\times
M\to M$}. Then \mbox{$\pr^{*}(\cP)\cong (K,\id\times \pi
,\bS^{n}\times P,\bS^{n}\times M)$}, where \mbox{$K$} acts trivially on the
first factor of \mbox{$\bS^{n}\times P$}. We have with respect to this action
\mbox{$C(\pr^{*}(P),K)^{K}\cong C (\bS^n\times P,K)^{K}$} and
\mbox{$C^{\infty}(\pr^{*}(P))^{K}\cong C^{\infty} (\bS^n\times P,K)^{K}$}.
The isomorphisms \mbox{$C(\bS^{n},G_{0})\cong C_{*}(\bS^{n},G_{0})\rtimes
G_{0}=C_{*}(\bS^{n},G)\rtimes G_{0}$}, where \mbox{$C_{*}(\bS^{n},G)$} denotes
the space of base-point-preserving maps from \mbox{$\bS^{n}$} to \mbox{$G$}, yield
\mbox{$\pi_{n}(G)=\pi_{0}(C_{*}(\bS^{n},G))= \pi_{0}(C(\bS^{n},G_{0}))$}
for any topological group \mbox{$G$}. We thus get a map
\begin{multline*}
\pi_{n}(C^{\infty}(P,K)^{K})=\pi_{0}(C_{*}(\bS^{n},C^{\infty}(P,K)^{K}))=\\
\pi_{0}(C(\bS^{n},C^{\infty}(P,K)^{K}_{\;\;0}))
\stackrel{\eta}{\to}\pi_{0}(C(\bS^{n},C(P,K)^{K}_{\;\; 0})),
\end{multline*}
where \mbox{$\eta$} is induced by the inclusion
\mbox{$C^{\infty}(P,K)^{K}\hookrightarrow C(P,K)^{K}$}.

If \mbox{$f\in C(\bS^{n}\times P,K)$} represents an element
\mbox{$[F]\in\pi_{0}(C(\bS^{n},C(P,K)^{K}_{\;\; 0}))$} (recall that we have
\mbox{$C(P,K)^{K}\cong G_{c,\cV }(\cP)\se \prod_{i=1}^{n}C(V_{i},K)$} and
\mbox{$C(\bS^{n},C(V_{i},K))\cong C(\bS^{n}\times V_{i},K)$}), then there
exists \mbox{$\wt{f}\in C^{\infty}(\bS^{n}\times P,K)^{K}$} which is
contained in the same connected component of \mbox{$C(\bS^{n}\times
P,K)^{K}$} as \mbox{$f$} (cf.\ \cite[Theorem 11]{approx}).  Since
\mbox{$\wt{f}$} is in particular smooth in the second argument, it follows
that \mbox{$\wt{f}$} represents an element \mbox{$\wt{F}\in
C(\bS^{n},C^{\infty}(P,K)^{K})$}.  Since the connected components and
the arc components of \mbox{$C(\bS^{n}\times P,K)^{K}$} coincide (since it is
a Lie group, cf. Remark \ref{rem:topologyOnContinuousGaugeGroup}),
there exists a path
\[
\tau :[0,1]\to C(\bS^{n}\times P,K)^{K}_{\;\;0}
\]
such that \mbox{$t\mapsto \tau (t)\cdot f$} is a path connecting \mbox{$f$} and
\mbox{$\wt{f}$}.  Since \mbox{$\bS^{n}$} is connected it follows that
\mbox{$C(\bS^{n}\times P,K)^{K}_{\;\;0}\cong C(\bS^{n},C(P,K)^{K})_{0}\se
C(\bS^{n},C(P,K)^{K}_{\;\;0})$}. Thus \mbox{$\tau$} represents a continuous path in
\mbox{$C(\bS^n,C(P,K)^{K}_{0}))$} connecting \mbox{$F$} and \mbox{$\wt{F}$} whence
\mbox{$[F]=[\wt{F}]\in \pi_{0}(C(\bS^{n},C(P,K)^{K}_{\;\;0}))$}.  That
\mbox{$\pi_{n} (\iota )$} is injective follows with Lemma
\ref{lem:approximationLemma2} as in \cite[Theorem A.3.7]{neeb03}.
\end{prf}

\section{The automorphism group as an infinite-dimensional Lie group}
\label{sect:theFullAutomorphismGroupAsInfiniteDimensionalLieGroup}

In this section we describe the Lie group structure on \mbox{$\Aut(\cP)$} for a
principal \mbox{$K$}-bundle over a compact manifold \mbox{$M$} \emph{without}
boundary, i.e., a \emph{closed} compact manifold. We will do this
using the extension of abstract groups
\begin{align}\label{eqn:extensionOfGauByDiff}
\Gau(\cP) \hookrightarrow \Aut(\cP)\xrightarrow{Q} \Diff(M)_{\cP},
\end{align}
where \mbox{$\Diff (M)_{\cP}$} is the image of the homomorphism \mbox{$Q\from
\Aut(\cP)\to\Diff(M)$}, \mbox{$F\mapsto F_{M}$} from Definition
\ref{def:bundleAutomorphismsAndGaugeTransformations}. 
More precisely, we will construct a Lie
group structure on \mbox{$\Aut(\cP)$} that turns
\eqref{eqn:extensionOfGauByDiff} into an extension of
Lie groups, i.e., into a locally trivial bundle.

We should advertise in advance that we shall not need regularity
assumptions on \mbox{$K$} in order to lift diffeomorphisms of \mbox{$M$} to bundle
automorphisms by lifting vector fields. However, we elaborate shortly
on the regularity of \mbox{$\Gau(\cP)$} and \mbox{$\Aut(\cP)$} in the end of the
section.
\vskip\baselineskip

We shall consider bundles over bases without boundary, i.e., our base
manifolds will always be closed compact manifolds. Throughout this
section we fix one particular given principal \mbox{$K$}-bundle \mbox{$\cP$} over a
closed compact manifold \mbox{$M$} and we furthermore assume that \mbox{$\cP$} has
the property SUB.

\begin{definition}(cf. \cite{neeb06nonAbelianExtensions})
If \mbox{$N$}, \mbox{$\wh{G}$} and \mbox{$G$} are Lie groups, then an extension of groups
\[
N\hookrightarrow \wh{G}\twoheadrightarrow G
\]
is called an \emph{extension of Lie groups} if \mbox{$N$} is a split Lie
subgroup of \mbox{$\wh{G}$}. That means that \mbox{$(N,q:\wh{G}\to G)$} is a smooth
principal \mbox{$N$}-bundle, where \mbox{$q:\wh{G}\to G\cong \wh{G}/N$} is the
induced quotient map.  We call two extensions \mbox{$N\hookrightarrow
\wh{G}_{1}\twoheadrightarrow G$} and \mbox{$N\hookrightarrow
\wh{G}_{2}\twoheadrightarrow G$} \emph{equivalent} if there exists a
morphism of Lie groups \mbox{$\psi \from \wh{G}_{1}\to \wh{G}_{2}$} such that
the diagram
\[
\begin{CD}
N@>>> \wh{G}_{1} @>>> G\\
@V\id_{N}VV @V{\psi}VV @V\id_{G}VV\\
N@>>> \wh{G}_{2} @>>> G
\end{CD}
\]
commutes.
\end{definition}

\begin{remark}\label{rem:choiceOfLocalTrivialisations} Unless stated
otherwise, for the rest of this section we choose and fix one
particular smooth closed trivialising system
\mbox{$\cl{\cV}=(\cl{V}_{i},\sigma_{i})_{\In}$} of \mbox{$\cP$} such that
\begin{itemize}
\item each \mbox{$\cl{V}_{i}$} is a compact manifold with corners diffeomorphic
	to \mbox{$[0,1]^{\dim(M)}$},
\item \mbox{$\cl{\cV}$} is a refinement of a smooth open trivialising system
	\mbox{${\mathcal{U}}=({U}_{i},\tau_{i})_{\In}$} and we have
	\mbox{$\cl{V}_{i}\se U_{i}$} and
	\mbox{$\sigma_{i}=\left.\tau_{i}\right|_{\cl{V}_{i}}$},
\item each \mbox{$\cl{U}_{i}$} is a compact manifold with corners diffeomorphic
	to \mbox{$[0,1]^{\dim(M)}$}
	and \mbox{$\tau_{i}$} extends to a smooth section
	\mbox{$\tau_{i}\from \cl{U}_{i}\to P$},
\item \mbox{$\cl{\mathcal{U}}\!=\!(\cl{U}_{i},\tau_{i})_{\In}$} is a refinement of a
	smooth open trivialising system \mbox{$\mathcal{U}'\!=\!(U'_{i},\tau_{j})_{\Jm}$}, 
\item the values of the transition functions
	\mbox{$k_{ij}\from U'_{i}\cap U'_{j}\to K$}
	of \mbox{$\mathcal{U}'$} are contained in open subsets \mbox{$W_{ij}$} of \mbox{$K$}, which
	are diffeomorphic to open zero neighbourhoods of \mbox{$\fk$},
\item \mbox{$\cP$} has the property SUB with respect to \mbox{$\cl{\cV}$} (and thus
with respect to \mbox{$\cl{\mathcal{U}}$} by Lemma
\ref{lem:propertySUBIsCompatibleWithRefinements}).
\end{itemize}

We choose \mbox{$\cl{\cV}$} by starting with an arbitrary smooth closed
trivialising system such that \mbox{$\cP$} has the property SUB with respect
to this system. Note that this exists because we assume throughout
this section that \mbox{$\cP$} has the property SUB.  Then Lemma
\ref{lem:forcingTransitionFunctionsIntoOpenCovers} implies that there
exists a refinement \mbox{$\mathcal{U}'=(U'_{j},\tau_{j})_{\Jm}$} such that
the transition functions \mbox{$k_{ij}\from U_{i}\cap U_{j}\to K$} take
values in open subsets \mbox{$W_{ij}$} of \mbox{$K$}, which are diffeomorphic to
open convex zero neighbourhoods of \mbox{$\fk$}.  Now each \mbox{$x\in M$} has
neighbourhoods \mbox{$V_{x}$} and \mbox{$U_{x}$} such that \mbox{$\cl{V}_{x}\se U_x$},
\mbox{$\cl{V}_{x}$} and \mbox{$\cl{U}_{x}$} are diffeomorphic to \mbox{$[0,1]^{\dim(M)}$}
and \mbox{$\cl{U}_{x}\se U_{j(x)}$} for some \mbox{$j(x)\in\{1,\dots ,m\}$}. Then
finitely many \mbox{$V_{x_{1}},\dots ,V_{x_{n}}$} cover \mbox{$M$} and so do
\mbox{$U_{x_{1}},\dots ,U_{x_{n}}$}. Furthermore, the sections \mbox{$\tau_{j}$}
restrict to smooth sections on \mbox{$V_{i}$}, \mbox{$\cl{V}_{i}$}, \mbox{$U_{i}$} and
\mbox{$\cl{U_{i}}$}.

This choice of \mbox{$\cl{\mathcal{U}}$} in turn implies that
\mbox{$\left.k_{ij}\right|_{\cl{U}_{i}\cap \cl{U}_{j}}$} arises as the
restriction of some smooth function on \mbox{$M$}.  In fact, if
\mbox{$\varphi_{ij}\from W_{ij}\to W'_{ij}\se \fk$} is a diffeomorphism onto
a convex zero neighbourhood and \mbox{$f_{ij}\in C^{\infty}(M,\R)$} is a
smooth function with \mbox{$\left.f_{ij}\right|_{\cl{U}_{i}\cap
\cl{U}_{j}}\equiv 1$} and \mbox{$\supp (f_{ij})\se U'_{i}\cap U'_{j}$}, then
\[
m\mapsto \left\{\begin{array}{ll}
\varphi_{ij}^{-1}(f_{ij}(m)\cdot \varphi_{ij}(k_{ij}(m)))
&	\text{ if }m\in U'_{i}\cap U'_{j}\\
\varphi_{ij}^{-1}(0)
&	\text{ if }m\notin U'_{i}\cap U'_{j}
\end{array}\right.
\]
is a smooth function, because each \mbox{$m\in \partial (U'_{i}\cap U'_{j})$}
has a neighbourhood on which \mbox{$f_{ij}$} vanishes, and this function
coincides with \mbox{$k_{ij}$} on \mbox{$\cl{U}_{i}\cap \cl{U}_{j}$}.

Similarly, let \mbox{$(\gamma_{1},\dots ,\gamma_{n})\in
G_{\cl{\mathcal{U}}}(\cP)\se \prod_{i=1}^{n}C^{\infty}(\cl{U}_{i},K)$}
be the local description of some element \mbox{$\gamma\in C^{\infty}(P,K)^{K}$}. We
will show that each \mbox{$\left.\gamma_{i}\right|_{\cl{V}_{i}}$} arises as
the restriction of a smooth map on \mbox{$M$}. In fact, take a diffeomorphism
\mbox{$\varphi_{i}\from \cl{U}_{i}\to [0,1]^{\dim(M)}$}. Then \mbox{$\cl{V}_{i}\se
U_{i}$} implies that we have \mbox{$\varphi_{i}(\cl{V}_{i})\se (0,1)^{\dim(M)}$} and thus
there exits an \mbox{$\varepsilon >0$} such that \mbox{$\varphi_{i} (\cl{V}_{i})\se
(\varepsilon ,1-\varepsilon)^{\dim (M)}$} for all \mbox{$i=1,\dots ,n$}. Now
let
\[
f\from [0,1]^{\dim(M)}\backslash (\varepsilon ,1-\varepsilon)^{\dim(M)}
\to [\varepsilon ,1-\varepsilon]^{\dim(M)} 
\]
be a map that restricts to the identity on \mbox{$\partial [\varepsilon
,1-\varepsilon]^{\dim(M)}$} and collapses \mbox{$\partial [0,1]^{\dim(M)}$} to a 
single point \mbox{$x_{0}$}. We then set
\[
\gamma'_{i}\from M\to K\quad m\mapsto \left\{\begin{array}{ll}
\gamma_{i}(m)
&	\text{ if }m\in \cl{U}_{i}\text{ and }
\varphi_{i} (m)\in [\varepsilon ,1-\varepsilon]^{\dim(M)}\\
\gamma_{i}(\varphi_{i}^{-1}(f(\varphi_{i} (m))))
&	\text{ if }m\in \cl{U}_{i}\text{ and }
\varphi_{i} (m)\notin (\varepsilon ,1-\varepsilon)^{\dim(M)}\\
\gamma_{i}(\varphi_{i}^{-1}(x_{0}))
&	\text{ if }m\notin U_{i},
\end{array} \right. 
\]
and \mbox{$\gamma '_{i}$} is well-defined and continuous, because
\mbox{$f(\varphi_{i} (m))=\varphi_{i} (m)$} if \mbox{$\varphi_{i}
(m)\in \partial [\varepsilon ,1-\varepsilon]^{\dim(M)}$} and
\mbox{$f(\varphi_{i} (m))=x_{0}$} if \mbox{$\varphi_{i} (m)\in\partial
[0,1]^{\dim(M)}$}. Since \mbox{$\gamma '_{i}$} coincides with \mbox{$\gamma_{i}$} on
the neighbourhood \mbox{$\varphi_{i}^{-1}((\varepsilon
,1-\varepsilon)^{\dim(M)})$}, it thus is smooth on this neighbourhood.
Now \cite[Corollary 12]{approx}, yields a
smooth map \mbox{$\wt{\gamma}_{i}$} on \mbox{$M$} with
\mbox{$\left.\gamma_{i}\right|_{\cl{V}_{i}}=
\left.\wt{\gamma}_{i}\right|_{\cl{V}_{i}}$}.
\end{remark}

We now give the description of a strategy for lifting special
diffeomorphisms to bundle automorphisms. This should motivate the
procedure of this section. \begin{remark}\label{rem:liftingDiffeomorphismsToBundleAutomorphisms}
Let \mbox{$U\se M$} be open and trivialising with section \mbox{$\sigma :U\to P$}
and corresponding \mbox{$k_{\sigma}:\pi^{-1}(U)\to K$}, given by \mbox{$\sigma (\pi
(p))\cdot k_{\sigma}(p)=p$}. If \mbox{$g \in\Diff(M)$} is such that
\mbox{$\supp(g)\se U$}, then we may define a smooth bundle automorphism
\mbox{$\wt{g}$} by
\[
\wt{g }(p)=\left\{
\begin{array}{ll}
\sigma \left(g\left(\pi  (p) \right) \right)\cdot k(p)
&\text{ if }p\in\pi^{-1}(U)\\
p&\text{else,}
\end{array}
\right.
\]
because each \mbox{$x\in \partial U$} has a neighbourhood on which \mbox{$\diffeo$}
is the identity. Furthermore, one easily verifies
\mbox{$Q(\wt{g})=\wt{g}_{M}=g$} and \mbox{$\wt{g^{-1}}=\wt{g}^{-1}$}, where
\mbox{$Q\from\Aut(\cP)\to \Diff(M)$} is the homomorphism from Definition
\ref{def:bundleAutomorphismsAndGaugeTransformations}.
\end{remark}

\begin{remark}\label{rem:chartsForDiffeomorphismGroups} Let \mbox{$M$} be a
closed compact manifold with a fixed Riemannian metric \mbox{$g$} and let
\mbox{$\pi \from TM\to M$} be its tangent bundle and \mbox{$\Exp \from TM\to M$} be
the exponential mapping of \mbox{$g$}. Then \mbox{$\pi \times \Exp\from TM\to
M\times M$}, \mbox{$X_{m}\mapsto (m,\Exp (X_{m}))$} restricts to a
diffeomorphism on an open neighbourhood \mbox{$U$} of the zero section in
\mbox{$TM$}.
We set \mbox{$O':=\{X\in\cV(M):X(M)\se U\}$} and define
\[
\varphi^{-1}:O'\to C^{\infty}(M,M),\quad \varphi^{-1} (X)(m)=\Exp(X(m))
\]
For the following, observe that \mbox{$\varphi^{-1}(X)(m)=m$} if and only if
\mbox{$X(m)=0_{m}$}.  After shrinking \mbox{$O'$} to a convex open neighbourhood in
the \mbox{$C^{1}$}-topology, one can also ensure that
\mbox{$\varphi^{-1}(X)\in \Diff(M)$} for all \mbox{$X\in O'$}. Since \mbox{$\pi
\times \Exp$} is bijective on \mbox{$U$}, \mbox{$\varphi^{-1}$} maps \mbox{$O'$} bijectively
to \mbox{$O:=\varphi^{-1}(O')\se \Diff(M)$} and thus endows \mbox{$O$} with a
smooth manifold structure. Furthermore, it can be shown that in view
of Proposition \ref{prop:localDescriptionsOfLieGroups}, this chart
actually defines a Lie group structure on \mbox{$\Diff(M)$}
(cf.\ \cite{leslie67}, \cite[Theorem 43.1]{krieglmichor97} or
\cite{gloeckner02patched}). It is even possible to put Lie group
structures on \mbox{$\Diff(M)$} in the case of non-compact manifolds,
possibly with corners \cite[Theorem 11.11]{michor80}, but we will not go
into this generality here.
\end{remark}

\begin{lemma}\label{lem:decompositionOfDiffeomorphisms} For the
open cover \mbox{$V_{1},\dots ,V_{n}$} of the closed compact manifold \mbox{$M$} and
the open identity neighbourhood \mbox{$O\se\Diff(M)$} from Remark
\ref{rem:chartsForDiffeomorphismGroups}, there exist smooth maps
\begin{align}\label{eqn:decompositionOfDiffeomorphisms}
s_{i}:O\to O\circ  O^{-1}
\end{align}
for \mbox{$1\leq i\leq n$} such that \mbox{$\supp(s_{i}(\diffeo))\se V_{i}$} and
\mbox{$s_{n}(g)\op{\circ }\dots \op{\circ } s_{1}(g)=g$}.
\end{lemma}

\begin{prf}(cf.\ \cite[Proposition 1]{hallerTeichmann04}) Let
\mbox{$f_{1},\ldots,f_{n}$} be a partition of unity subordinated to the open
cover \mbox{$V_{1},\ldots,V_{n}$} and let \mbox{$\varphi\from O\to \varphi (O)\se
\cV (M)$} be the chart of \mbox{$\Diff(M)$} form Remark
\ref{rem:chartsForDiffeomorphismGroups}. In particular,
\mbox{$\varphi^{-1}(X)(m)=m$} if \mbox{$X(m)=0_{m}$}. Since \mbox{$\varphi (O)$} is
convex, we may define \mbox{$s_{i}\from O\to O\op{\circ}O^{-1}$},
\[
s_{i}(\diffeo)= 
\varphi^{-1}\big((f_{n}+\ldots+f_{i})\cdot\varphi (g) \big)
\op{\circ}
\big(\varphi^{-1}\big((f_{n}+\ldots+f_{i+1})\cdot \varphi(g) \big)\big)^{-1}
\]
if \mbox{$i<n$} and \mbox{$s_{n}(\diffeo)=\varphi^{-1} (f_{n}\cdot \varphi (\diffeo
))$}, which are smooth since they are given by a push-forward of the
smooth map \mbox{$\R \times TM\to TM$} \mbox{$(\lambda ,X_{m})\mapsto \lambda \cdot
X_{m}$}. Furthermore, if \mbox{$f_{i}(x)=0$}, then the left and the right
factor annihilate each other and thus \mbox{$\supp(s_{i}(\diffeo ))\se
V_{i}$}.
\end{prf}

The preceding lemma enables us now to lift elements
of \mbox{$O\se \Diff(M)$} to elements of \mbox{$\Aut(\cP)$}.

\begin{definition}\label{def:sectionFromDiffeomorphismsToBundleAutomorphisms}
If \mbox{$O\se\Diff(M)$} is the open identity neighbourhood from Remark
\ref{rem:chartsForDiffeomorphismGroups} and \mbox{$s_{i}:O\to O\op{\circ
}O^{-1}$} are the smooth mappings from
Lemma \ref{lem:decompositionOfDiffeomorphisms}, then we define 
\begin{align}\label{eqn:sectionFromDiffeomorphismIntoBundleAutoporphisms}
S\from O\to \Aut(\cP),\quad \diffeo \mapsto 
S(\diffeo ):=\wt{\diffeo_{n}}\op{\circ}\dots \op{\circ}\wt{\diffeo_{1\,}},
\end{align}
where \mbox{$\wt{\diffeo_{i}}$} is the bundle
automorphism of \mbox{$\cP$} from Remark
\ref{rem:liftingDiffeomorphismsToBundleAutomorphisms}.  This defines a
local section for the homomorphism \mbox{$Q\from \Aut(\cP)\to \Diff(M)$},
\mbox{$F\mapsto F_{M}$} from Definition
\ref{def:bundleAutomorphismsAndGaugeTransformations}.
\end{definition}

We shall frequently need an explicit description of \mbox{$S(\diffeo)$} in
terms of local trivialisations, i.e., how \mbox{$S(\diffeo)(\sigma_{i}(x))$}
can be expressed in terms of \mbox{$\diffeo_{j}$}, \mbox{$\sigma_{j}$} and
\mbox{$k_{jj'}$}. 

\begin{remark}\label{rem:valuesOfTheLiftedDiffeomorphismInTermsOfLocalData}
Let \mbox{$x\in V_{i}\se M$} be such that \mbox{$x\notin V_{j}$} for \mbox{$j<i$} and
\mbox{$g_{i}(x)\notin V_{j}$} for \mbox{$j> i$}. Then \mbox{$\diffeo_{j}(x)=x$} for all
\mbox{$j<i$}, \mbox{$\diffeo_{j} (\diffeo_{i}(x))=\diffeo_{i}(x)$} for all \mbox{$j>i$} and
thus \mbox{$S(\diffeo
)(\sigma_{i}(x))=\sigma_{i}(\diffeo_{i}(x))=\sigma_{i}(\diffeo (x))$}.

In general, things are more
complicated. The first \mbox{$\wt{g_{j_{1}}}$} in
\eqref{eqn:sectionFromDiffeomorphismIntoBundleAutoporphisms} that
could move \mbox{$\sigma_{i}(x)$} is the one for the minimal \mbox{$j_{1}$} such
that \mbox{$x\in \cl{V}_{j_{1}}$}. We then have 
\[
 \wt{\diffeo_{j_{1}}}(\sigma_{i}(x))
=\wt{\diffeo_{j_{1}}}(\sigma_{j_{1}}(x))\cdot k_{j_{1}i}(x)
=\sigma_{j_{1}}({\diffeo_{j_{1}}}(x))\cdot k_{j_{1}i}(x).
\]
The next \mbox{$\wt{\diffeo_{j_{2}}}$} in
\eqref{eqn:sectionFromDiffeomorphismIntoBundleAutoporphisms} that
could move \mbox{$\wt{\diffeo_{j_{1}}}(\sigma_{i}(x))$} in turn is the one
for the minimal \mbox{$j_{2}>j_{1}$} such that \mbox{$\diffeo_{j_{1}}(x)\in
\cl{V}_{j_{2}}$}, and we then have
\[
 \wt{\diffeo_{j_{2}}}(\wt{\diffeo_{j_{1}}}(\sigma_{i}(x)))
=\sigma_{j_{2}}(\diffeo_{j_{2}}\op{\circ }{\diffeo_{j_{1}}}(x))
\cdot k_{j_{2}j_{1}}(\diffeo_{j_{1}}(x))\cdot k_{j_{1}i}(x).
\]
We eventually get
\begin{align}\label{eqn:valuesOfTheLiftedDiffeomorphismsInTermsOflocalData}
S(\diffeo )(\sigma_{i}(x)) =\sigma_{j_{\ell}}(\diffeo (x)) \cdot
k_{j_{\ell}j_{\ell-1}}(\diffeo_{j_{\ell-1}}\op{\circ}\dots \op{\circ}
\diffeo_{j_{1}}(x))\cdot\ldots \cdot  k_{j_{1}i}(x),
\end{align}
where \mbox{$\{j_{1},\dots ,j_{\ell}\}\se\{1,\dots ,n\}$} is maximal such
that 
\[
\diffeo_{j_{p-1}}\op{\circ}\ldots\op{\circ}\diffeo_{i_{1}}(x)\in
U_{j_{p}}\cap U_{j_{p-1}}
\;\text{ for }\; 2\leq p\leq \ell\;\text{ and }\;j_{1}<\ldots
<j_{p}.
\]

Note that we cannot
write down such a formula using all \mbox{$j\in\{1,\dots ,n\}$}, because
the corresponding \mbox{$k_{jj'}$} and \mbox{$\sigma_{j}$} would not be defined properly.

Of course, \mbox{$\diffeo$} and \mbox{$x$} influence the choice of \mbox{$j_{1},\dots
,j_{\ell}$}, but there exist open neighbourhoods \mbox{$O_{g}$} of \mbox{$\diffeo$}
and \mbox{$U_{x}$} of \mbox{$x$} such that we may use
\eqref{eqn:valuesOfTheLiftedDiffeomorphismsInTermsOflocalData} as a
formula for all \mbox{$\diffeo '\in O_{\diffeo}$} and \mbox{$x'\in U_{x}$}.  In
fact, the action \mbox{$\Diff(M)\times M\to M$}, \mbox{$\diffeo.m=\diffeo (m)$} is
smooth by \cite[Proposition 7.2]{gloeckner02patched}, and thus in
particular continuous. If
\begin{align}
\diffeo_{j_{p}}\op{\circ}\ldots\op{\circ}\diffeo_{j_{1}}(x)\notin
\cl{V}_{j} \;&\text{ for }\;
2\leq p\leq \ell\;\text{ and }\;j\notin \{j_{1},\dots ,j_{p}\}
\label{eqn:valuesOfTheLiftedDiffeomorphismsInTermsOflocalData1}\\
\diffeo_{j_{p}}\op{\circ}\ldots\op{\circ}\diffeo_{j_{1}}(x)\in
U_{j_{p}}\cap U_{j_{p-1}}
\;&\text{ for }\; 2\leq p\leq \ell\;\text{ and }\;j_{1}<\ldots
<j_{p}
\label{eqn:valuesOfTheLiftedDiffeomorphismsInTermsOflocalData2}
\end{align}
then this is also true for \mbox{$\diffeo '$} and \mbox{$x'$} in some open
neighbourhood of \mbox{$\diffeo$} and \mbox{$x$}. This yields finitely many
open neighbourhoods of \mbox{$\diffeo$} and \mbox{$x$} and we define their
intersections to be \mbox{$O_{\diffeo }$} and \mbox{$U_{x}$}. Then
\eqref{eqn:valuesOfTheLiftedDiffeomorphismsInTermsOflocalData} still
holds for \mbox{$\diffeo '\in O_{\diffeo }$} and \mbox{$x'\in U_{x}$}, because
\eqref{eqn:valuesOfTheLiftedDiffeomorphismsInTermsOflocalData1}
implies
\mbox{$\diffeo_{j}(\diffeo_{j_{p}}\op{\circ}\ldots\op{\circ}\diffeo_{j_{1}}(x))
=\diffeo_{j_{p}}\op{\circ}\ldots\op{\circ}\diffeo_{j_{1}}(x)$} and
\eqref{eqn:valuesOfTheLiftedDiffeomorphismsInTermsOflocalData2}
implies that \mbox{$k_{j_{p}j_{p-1}}$} is defined and satisfies the cocycle
condition.
\end{remark}

In order to determine a Lie group structure on \mbox{$\Aut(\cP)$}, the map
\mbox{$S\from O\to \Aut(\cP)$} has to satisfy certain smoothness
properties, which will be ensured by the subsequent lemmas.

\begin{remark}\label{rem:actionsOfTheAutomorphismOnMappingsInterchange}
If we identify the normal subgroup \mbox{$\Gau(\cP)\unlhd \Aut(\cP)$} with
\mbox{$C^{\infty }(P,K)^{K}$} via
\[
C^{\infty}(P,K)^{K}\to \Gau(\cP),\quad
\gamma \mapsto F_{\gamma }
\]
with \mbox{$F_{\gamma}(p)=p\cdot \gamma (p)$}, then the conjugation action
\mbox{$c:\Aut(\cP)\times \Gau(\cP)\to\Gau(\cP)$}, given by
\mbox{$c_{F}(F_{\gamma })=F\op{\circ}F_{\gamma}\op{\circ}F^{-1}$} changes into
\[
c:\Aut(\cP)\times C^{\infty}(P,K)^{K}\to C^{\infty}(P,K)^{K},\quad
(F,\gamma)\mapsto \gamma \op{\circ}F^{-1}.
\]
In fact, this follows from
\[
(F\circ F_{\gamma}\circ F^{-1})(p)
=F\big(F^{-1}(p)\cdot \gamma (F^{-1}(p))\big)
=p\cdot \gamma (F^{-1}(p))=F_{(\gamma \circ F^{-1})}(p).
\]
\end{remark}

In the following remarks and lemmas we show the smoothness of the
maps \mbox{$T$}, \mbox{$\omega$} and \mbox{$\omega_{\diffeo}$}, mentioned before.

\begin{lemma}\label{lem:localActionOftheDiffeomorphismGroupOnTheGaugeAlgebra}
Let \mbox{$O\se \Diff(M)$} be the open identity neighbourhood from
Remark \ref{rem:chartsForDiffeomorphismGroups} and let \mbox{$S:O\to
\Aut(\cP)$} be the map from Definition
\ref{def:sectionFromDiffeomorphismsToBundleAutomorphisms}. Then we
have that for each \mbox{$F\in \Aut(\cP)$} the map
\mbox{$C^{\infty}(P,\fk )^{K}\to C^{\infty}(P,\fk )^{K}$}, \mbox{$\eta
\mapsto \eta \op{\circ}F^{-1}$} is an automorphism of
\mbox{$C^{\infty}(P,\fk)^{K}$} and the map
\[
t\from C^{\infty}(P,\fk)^{K}\times O\to C^{\infty}(P,\fk)^{K},\quad
(\eta ,\diffeo )\mapsto \eta \op{\circ}S(\diffeo )^{-1}
\]
is smooth.
\end{lemma}

\begin{prf}
That \mbox{$\eta \mapsto \eta \op{\circ}F^{-1}$} is an element of
\mbox{$\Aut(C^{\infty}(P,\fk)^{K})$} follows immediately from the (pointwise)
definition of the bracket on \mbox{$C^{\infty}(P,\fk)^{K}$}. We shall use the
previous established isomorphisms \mbox{$C^{\infty}(P,\fk)^{K}\cong
\fg_{\mathcal{U}'}(\cP)\cong\fg_{\mathcal{U}}(\cP)\cong\fg_{\cV}(\cP)$}
from Proposition \ref{prop:isomorphismOfTheGaugeAlgebra} and reduce
the smoothness of \mbox{$t$} to the smoothness of
\[
C^{\infty}(M,\fk)\times \Diff(M)\to C^{\infty}(M,\fk),\quad
(\eta ,\diffeo)\mapsto \eta \op{\circ }\diffeo^{-1}
\]
from \cite[Proposition 6]{gloeckner02patched} and to
the action of \mbox{$\diffeo_{i}^{-1}$} on \mbox{$C^{\infty}(\cl{V}_{i},\fk)$},
because we have no description of what \mbox{$\diffeo_{i}^{-1}$} does with
\mbox{$U_{j}$} for \mbox{$j\neq i$}. It clearly suffices to show that the map
\[
t_{i}:C^{\infty}(P,\fk)^{K}\times \Diff(M)\to
C^{\infty}(P,\fk)^{K}\times \Diff (M),\quad (\eta ,\diffeo )\mapsto  
(\eta\op{\circ}\wt{g_{i}}^{-1},\diffeo )
\]
is smooth for each \mbox{$1\leq i\leq n$}, because then \mbox{$t=\pr_{1}\op{\circ
}t_{n}\op{\circ}\ldots\op{\circ}t_{1}$} is smooth. This in turn follows
from the smoothness of
\begin{align}\label{eqn:localActionOnTheCurrentAlgebra}
C^{\infty}(U'_{i},\fk)\times \Diff(M)\to C^{\infty}({U}_{i},\fk),\quad
(\eta ,\diffeo)\mapsto \eta \op{\circ }
\left.\diffeo_{i}^{-1}\right|_{U_{i}},
\end{align}
because this is the local description of \mbox{$t_{i}$}. In fact, for each
\mbox{$j\neq i$} there exists an open subset \mbox{$V'_{j}$} with \mbox{$U_{j}\backslash
U_{i}\se V'_{j}\se U_{j}\backslash V_{i}$}, because \mbox{$\cl{V}_{i}\se
U_{i}$} and \mbox{$U_{j}$} is diffeomorphic to \mbox{$(0,1)^{\dim(M)}$}.
Furthermore, we set \mbox{$V'_{i}:=U_{i}$}. Then \mbox{$(V'_{1},\dots ,V'_{n})$} is
an open cover of \mbox{$M$}, leading to a refinement \mbox{$\cV'$} of the
trivialising system \mbox{$\mathcal{U}'$} and we have
\[
t_{i}:\fg_{\mathcal{U}'}(\cP)\times O\to \fg_{\cV '}(\cP),\quad
((\eta_{1},\dots ,\eta_{n}),\diffeo )\mapsto 
(\left.\eta_{1}\right|_{{V'_{1}}},\dots ,
 \left.\eta_{i}\op{\circ}g_{i}^{-1}\right|_{{V'_{i}}},\dots ,
 \left.\eta_{n}\right|_{{V'_{n}}})
\]
because \mbox{$\supp (\diffeo_{i})\se V_{i}$} and \mbox{$V'_{j}\cap
V_{i}=\emptyset$} if \mbox{$j\neq i$}.  To show that
\eqref{eqn:localActionOnTheCurrentAlgebra} is smooth, choose some
\mbox{$f_{i}\in C^{\infty}(M,\R)$} with \mbox{$\left.f_{i}\right|_{U_{i}}\equiv 1$}
and \mbox{$\supp (f_{i})\se U'_{i}$}. Then
\[
h_{i}:C^{\infty}(U'_{i},\fk)\to C^{\infty}(M,\fk ),\quad
\eta \mapsto \left( m\mapsto \left\{\begin{array}{ll}
f_{i}(m)\cdot \eta (m) & \text{ if } m\in U'_{i}\\
0 & \text{ if }m\notin U'_{i}
\end{array} \right.\right)
\]
is smooth by Corollary \ref{cor:gluingForVectorBundles_OPEN_Version},
because \mbox{$\eta \mapsto \left.f_{i}\right|_{U'_{i}}\cdot\eta $} is
linear, continuous and thus smooth. Now we have \mbox{$\supp(\diffeo_{i})\se
V_{i}\se U_{i}$} and thus
\mbox{$\left.h_{i}(\eta)\op{\circ}g_{i}^{-1}\right|_{U_{i}}
=\left.\eta\op{\circ}g_{i}^{-1}\right|_{U_{i}}$}
depends smoothly on \mbox{$\diffeo$} and \mbox{$\eta$} by
Corollary
\ref{cor:restrictionMapIsSmooth}.
\end{prf}

The following proofs share a common idea. We will always have to show
that certain mappings with values in \mbox{$C^{\infty}(P,K)^{K}$} are smooth.
This can be established by showing that their compositions with the
pull-back \mbox{$(\sigma_{i})^{*}$} of a section \mbox{$\sigma_{i}\from \cl{V}_{i}\to
P$} (then with values in \mbox{$C^{\infty}(\cl{V}_{i},K)$}) are smooth for all
\mbox{$1\leq i \leq n$}.  As described in Remark
\ref{rem:valuesOfTheLiftedDiffeomorphismInTermsOfLocalData}, it will
not be possible to write down explicit formulas for these mappings in
terms of the transition functions \mbox{$k_{ij}$} for all \mbox{$x\in\cl{V}_{i}$}
simultaneously, but we will be able to do so on some open
neighbourhood \mbox{$U_{x}$} of \mbox{$x$}.  For different \mbox{$x_{1}$} and \mbox{$x_{2}$} these
formulas will define the same mapping on \mbox{${U}_{x_{1}}\cap
{U}_{x_{2}}$}, because there they define
\mbox{$(\sigma_{i}^{*}(S(\diffeo)))=S(\diffeo)\op{\circ }\sigma_{i}$}. By
restriction and gluing we will thus be able to reconstruct the
original mappings and then see that they depend smoothly on their
arguments.

\begin{lemma}
\label{lem:orbitMapFromDiffeomorhismsToBundleAutomorphismsIsSmooth}
If \mbox{$O\se \Diff(M)$} is the open identity neighbourhood from
Remark \ref{rem:chartsForDiffeomorphismGroups} and if \mbox{$S:O\to \Aut(\cP)$}
is the map from Definition
\ref{def:sectionFromDiffeomorphismsToBundleAutomorphisms},
then for each \mbox{$\gamma \in C^{\infty}(P,K)^{K}$} the map
\[
O\ni \diffeo \mapsto \gamma \op{\circ}S(\diffeo)^{-1}\in C^{\infty}(P,K)^{K}
\]
is smooth.
\end{lemma}

\begin{prf}
It suffices to show that \mbox{$\gamma \op{\circ}S(\diffeo)^{-1}\op{\circ
}\left.\sigma_{i}\right|_{\cl{V}_{i}}$} depends smoothly on \mbox{$\diffeo$}
for \mbox{$1\leq i\leq n$}.  Let \mbox{$(\gamma_{1},\dots ,\gamma_{n})\in
G_{\mathcal{U}}(\cP)\se \prod_{i=1}^{n}C^{\infty}(\cl{U}_{i},K)$} be
the local description of \mbox{$\gamma$}.  Fix \mbox{$\diffeo \in O$} and \mbox{$x\in
\cl{V}_{i}$}. Then Remark
\ref{rem:valuesOfTheLiftedDiffeomorphismInTermsOfLocalData} yields
open neighbourhoods \mbox{$O_{\diffeo}$} of \mbox{$\diffeo$} and \mbox{$U_{x}$} of \mbox{$x$}
(w.l.o.g. such that \mbox{$\cl{U}_{x}\se \cl{V}_{i}$} is a manifold with
corners) such that
\begin{multline*}
 \gamma (S(\diffeo')^{-1}(\sigma_{i}(x')))
= \gamma \big(\sigma_{j_{\ell}}(\diffeo' (x'))
\underbrace{\cdot k_{j_{\ell}j_{\ell-1}}(\diffeo'_{j_{\ell-1}}\op{\circ}\dots
\op{\circ} \diffeo'_{j_{1}}(x'))\cdot\ldots \cdot
k_{j_{1}i}(x')}_{:=\kappa_{x,\diffeo '}(x')} \big)\\
= \kappa_{x,\diffeo'}(x')^{-1}\cdot\gamma\big(\sigma_{j_{\ell}}(\diffeo'(x'))
  \cdot
  \kappa_{x,\diffeo'}(x') 
= \underbrace{
  \kappa_{x,\diffeo'}(x')^{-1}\cdot\gamma_{j_{\ell}}(\diffeo'(x'))\cdot
  \kappa_{x,\diffeo'}(x')}_{:=\theta_{x,\diffeo'}(x')} 
\end{multline*}
for all \mbox{$\diffeo '\in O_{\diffeo}$} and \mbox{$x'\in \cl{U}_{x}$}. Since we
will not vary \mbox{$i$} and \mbox{$\diffeo$} in the sequel, we suppressed the
dependence of \mbox{$\kappa_{x,\diffeo '}(x')$} and \mbox{$\theta_{x,\diffeo
'}(x')$} on \mbox{$i$} and \mbox{$\diffeo$}. Note that each \mbox{$k_{jj'}$} and
\mbox{$\gamma_{i}$} can be assumed to be defined on \mbox{$M$} (cf.\ Remark
\ref{rem:choiceOfLocalTrivialisations}). Thus, for fixed \mbox{$x$}, the
formula for \mbox{$\theta _{x,\diffeo '}$} defines a smooth function on \mbox{$M$}
that depends smoothly on \mbox{$g'$}, because the action of \mbox{$\Diff(M)$} on
\mbox{$C^{\infty}(M,K)$} is smooth (cf.\ \cite[Proposition
10.3]{gloeckner02patched}).

Furthermore, \mbox{$\theta_{x_{1},\diffeo'}$} and \mbox{$\theta_{x_{2},\diffeo '}$}
coincide on \mbox{$\cl{U}_{x_{1}}\cap \cl{U}_{x_{2}}$}, because
both define \mbox{$\gamma \op{\circ}S(\diffeo ')^{-1}\op{\circ }\sigma_{i}$} there. Now
finitely many \mbox{$U_{x_{1}},\dots ,U_{x_{m}}$} cover \mbox{$\cl{V}_{i}$}, and 
since the gluing and restriction maps from Lemma
\ref{lem:restrictionMapForCurrentGroupIsSmooth} and Proposition
\ref{prop:gluingLemmaForCurrentGroup} are smooth we
have that
\[
\gamma \op{\circ}S(\diffeo ')^{-1}\op{\circ}\sigma_{i}=\glue
(\left.\theta_{x_{1},\diffeo '}\right|_{\cl{U}_{x_{1}}}
 ,\dots ,
 \left.\theta_{x_{m},\diffeo '}\right|_{\cl{U}_{x_{m}}})
\]
depends smoothly on \mbox{$\diffeo '$}.
\end{prf}

\begin{lemma}\label{lem:localActionOfTheDiffeomorphismGroupOnTheGaugeGroup}
Let \mbox{$O\se \Diff(M)$} be the open identity neighbourhood from Remark
\ref{rem:chartsForDiffeomorphismGroups} and let \mbox{$S:O\to \Aut(\cP)$} be the
map from Definition
\ref{def:sectionFromDiffeomorphismsToBundleAutomorphisms}. Then we have that for each
\mbox{$F\in \Aut(\cP)$} the map \mbox{$c_{F}\from C^{\infty}(P,K)^{K}\to
C^{\infty}(P,K)^{K}$}, \mbox{$\gamma \mapsto \gamma \op{\circ}F^{-1}$} is an
automorphism of \mbox{$C^{\infty}(P,K)^{K}$} and the map
\begin{align}
\label{eqn:definitionOfTheOuterActionOfTheDiffeomorphismsOnTheGaugeGroup}
T:C^{\infty}(P,K)^{K}\times O\to C^{\infty }(P,K)^{K},\quad
(\gamma ,\diffeo)\mapsto \gamma \op{\circ}S(g)^{-1}
\end{align}
is smooth.
\end{lemma}

\begin{prf}
Since \mbox{$\gamma \mapsto \gamma \op{\circ}F^{-1}$} is a group
homomorphism, it suffices to show that it is smooth on a unit
neighbourhood. Because the charts on \mbox{$C^{\infty}(P,K)^{K}$} are
constructed by push-forwards (cf.\ Proposition
\ref{prop:gaugeGroupInLocalCoordinatesIsLieGroup}) this follows
immediately from the fact that the corresponding automorphism of
\mbox{$C^{\infty}(P,\fk)^{K}$}, given by \mbox{$\eta \mapsto \eta
\op{\circ}F^{-1}$}, is continuous and thus smooth. For the same reason,
Lemma \ref{lem:localActionOftheDiffeomorphismGroupOnTheGaugeAlgebra}
implies that there exists a unit neighbourhood \mbox{$U\se
C^{\infty}(P,K)^{K}$} such that
\[
U\times O\to C^{\infty}(P,K)^{K},\quad (\gamma ,\diffeo)\mapsto \gamma
\op{\circ}S(\diffeo )^{-1}
\] 
is smooth.

Now for each \mbox{$\gamma_{0} \in C^{\infty}(P,K)^{K}$} there exists an open
neighbourhood \mbox{$U_{\gamma_{0} }$} with \mbox{$\gamma_{0}^{-1}\cdot
U_{\gamma_{0}}\se U$}. Hence
\[
\gamma \op{\circ}S(\diffeo )^{-1}=(\gamma_{0}\cdot
\gamma_{0}^{-1}\cdot \gamma )\op{\circ}S(\diffeo)^{-1}=
\big(\gamma_{0}\op{\circ}S(\diffeo)^{-1}\big)\cdot\big(
(\gamma_{0}^{-1}\cdot \gamma)\op{\circ}S(\diffeo)^{-1}\big),
\]
and the first factor depends smoothly on \mbox{$\diffeo$} due to Lemma
\ref{lem:orbitMapFromDiffeomorhismsToBundleAutomorphismsIsSmooth}, and
the second factor depends smoothly on \mbox{$\gamma$} and \mbox{$\diffeo$}, because
\mbox{$\gamma_{0}^{-1}\cdot \gamma \in U$}.
\end{prf}

\begin{lemma}\label{lem:smoothCocycleOnAutomorhismGroup}
If \mbox{$O\se \Diff(M)$} is the open identity neighbourhood from
Remark \ref{rem:chartsForDiffeomorphismGroups} and if \mbox{$S:O\to\Aut(\cP)$}
is the map from Definition
\ref{def:sectionFromDiffeomorphismsToBundleAutomorphisms}, then 
\begin{align}
\label{eqn:definitionOfTheCocycleFromDiffeomorphismsToTheGaugeGroup}
\omega \from O\times O\to C^{\infty}(P,K)^{K},\quad 
(\diffeo ,\diffeo ')\mapsto
S (\diffeo )\op{\circ } S (\diffeo ')\op{\circ } S(\diffeo \op{\circ }
\diffeo')^{-1}
\end{align}
is smooth. Furthermore, if \mbox{$Q\from \Aut(\cP)\to \Diff(M)$}, \mbox{$F\mapsto
F_{M}$} is the homomorphism from Definition
\ref{def:bundleAutomorphismsAndGaugeTransformations} then for each
\mbox{$\diffeo \in Q(\Diff(M))$} there exists an open identity neighbourhood
\mbox{$O_{\diffeo}\se O$} such that
\begin{align}\label{eqn:conjugationOfTheCocycle}
\omega_{g}:O_{\diffeo }\to C^{\infty}(P,K)^{K},\quad  
\diffeo '\mapsto F\op{\circ } S(\diffeo ')\op{\circ } F^{-1}
\op{\circ } S(\diffeo \op{\circ } \diffeo' \op{\circ }\diffeo^{-1})^{-1}
\end{align}
is smooth for any \mbox{$F\in \Aut(\cP)$} with \mbox{$F_{M}=g$}.
\end{lemma}

\begin{prf}
First observe that \mbox{$\omega (\diffeo ,\diffeo ')$}
actually is an element of \mbox{$C^{\infty}(P,K)^{K}\cong
\Gau(\cP)=\ker(Q)$}, because \mbox{$Q$} is a homomorphism of groups, \mbox{$S$} is a
section of \mbox{$Q$} and thus
\[
Q(\omega (\diffeo ,\diffeo '))=Q(S(\diffeo))\op{\circ } Q(S(\diffeo
'))\op{\circ } Q(S(\diffeo \op{\circ}\diffeo '))^{-1}=\id_{M}.
\]

To show that \mbox{$\omega$} is smooth, we derive an explicit formula for
\mbox{$\omega(\diffeo, \diffeo ')\op{\circ}\sigma_{i}\in
C^{\infty}(\cl{V}_{i},K)$} that depends smoothly on \mbox{$\diffeo$} and
\mbox{$\diffeo '$}.

Denote \mbox{$\wh{\diffeo }:=\diffeo \op{\circ}\diffeo '$} for \mbox{$\diffeo
,\diffeo '\in O$} and fix \mbox{$\diffeo ,\diffeo '\in O$}, \mbox{$x\in
\cl{V}_{i}$}. Proceeding as in Remark
\ref{rem:valuesOfTheLiftedDiffeomorphismInTermsOfLocalData}, we find
\mbox{$i_{1},\dots ,i_{\ell}$} such that
\[
S(\wh{\diffeo})^{-1}(\sigma_{i_{\ell}}(x))=\sigma_{\ell}(\wh{\diffeo
}^{-1}(x))\cdot
k_{i_{\ell}i_{\ell-1}}((\wh{\diffeo}_{i_{\ell-1}})^{-1}\op{\circ}\ldots
\op{\circ}\left(\wh{\diffeo}_{i_{1}}\right)^{-1}(x))\cdot \ldots\cdot
k_{i_{1}i}(x).
\]
Accordingly we find \mbox{$i'_{\ell'},\dots ,i'_{1}$} for \mbox{$S(\diffeo')$} and
\mbox{$i''_{\ell''},\dots ,i''_{1}$} for \mbox{$S(\diffeo)$}. We get as in Remark
\ref{rem:valuesOfTheLiftedDiffeomorphismInTermsOfLocalData} open
neighbourhoods \mbox{$O_{\diffeo },O_{\diffeo '}$} of \mbox{$\diffeo ,\diffeo '$}
and \mbox{$U_{x}$} of \mbox{$x$} (w.l.o.g. such that \mbox{$\cl{U}_{x}\se \cl{V}_{i}$} is a
manifold with corners) such that for \mbox{$h\in O_{\diffeo}$}, \mbox{$h\in
O_{\diffeo '}$} and \mbox{$x'\in \cl{U}_{x}$} we have \mbox{$ S(h )\cdot S(h ')\cdot
S(h \cdot h')^{-1}(\sigma_{i}(x')) = $}
\begin{alignat*}{3}
\sigma_{i}(x') \cdot \Big[ & k_{i\,i''_{\ell''}}(x') &&&&\\
  \cdot 
& k_{i''_{\ell''}i''_{\ell''-1}}\big(h_{i''_{\ell''-1}}\op{\circ}\ldots
  \op{\circ }h_{i''_{1}}\op{\circ } h^{-1}(x')\big)
  \cdot\ldots & \,\cdot
& k_{i''_{1}i'_{\ell'}}(h^{-1}(x')) \\
  \cdot
& k_{i'_{\ell'}i'_{\ell'-1}}\big(h'_{i'_{\ell'-1}}\op{\circ }\ldots
  \op{\circ} h'_{i'_{1}}\op{\circ } \wh{h}^{-1}(x')\big)
  \cdot \ldots & \cdot 
& k_{i'_{1}i_{\ell}}(\wh{h}^{-1}(x')) \\
  \cdot 
& k_{i_{\ell}i_{\ell-1}}\big((\wh{h}_{i_{\ell-1}})^{-1}\op{\circ}\ldots
  \op{\circ }\left.(\wh{h}_{i_{1}}\right.)^{-1}(x')\big) 
  \cdot\ldots & \cdot
& k_{i_{1}i}(x')\Big].
\end{alignat*}
Denote by \mbox{$\kappa_{x,h,h'}(x')\in K$} the element in brackets on the
right hand side, and note that it defines \mbox{$\omega
(h,h')\op{\circ}\sigma_{i}(x')$} by Remark
\ref{rem:gaugeGroupIsIsomorphicToEquivariantMappings}. Since we will
not vary \mbox{$\diffeo$} and \mbox{$\diffeo '$} in the sequel we suppressed
the dependence of \mbox{$\kappa_{x,h,h'}(x')$} on them.

Now each \mbox{$k_{ij}$} can be assumed to be defined on \mbox{$M$} (cf.\ Remark
\ref{rem:choiceOfLocalTrivialisations}). Thus, for fixed \mbox{$x$}, the
formula for \mbox{$\kappa_{x,h,h'}$} defines a smooth function on \mbox{$M$} that
depends smoothly on \mbox{$h$} and \mbox{$h'$}, because the action of \mbox{$\Diff(M)$} on
\mbox{$C^{\infty}(M,K)$} is smooth (cf.\ \cite[Proposition
10.3]{gloeckner02patched}).  Furthermore, \mbox{$\kappa_{x_{1},h,h'}$}
coincides with \mbox{$\kappa_{x_{2},h,h'}$} on \mbox{$\cl{U}_{x_{1}}\cap
\cl{U}_{x_{2}}$}, because
\[
\sigma_{i}(x')\cdot \kappa_{x_{1},h,h'}(x')=S(h )\op{\circ } S(h
')\op{\circ } S(h \op{\circ } h')^{-1}(\sigma_{i}(x')) =\sigma_{i}(x')\cdot
\kappa_{x_{2},h,h'}(x')
\]
for \mbox{$x'\in \cl{U}_{x_{1}}\cap \cl{U}_{x_{2}}$}.
Now finitely many \mbox{$U_{x_{1}},\dots ,U_{x_{m}}$} cover \mbox{$\cl{V}_{i}$} and we
thus see that
\[
\omega (h,h')\op{\circ }\sigma_{i}=
\glue(
\left.\kappa_{x_{1},h,h'}\right|_{\cl{U}_{x_{1}}}
,\dots ,
\left.\kappa_{x_{m},h,h'}\right|_{\cl{U}_{x_{m}}}
)
\]
depends smoothly on \mbox{$h$} and \mbox{$h'$}.

We derive an explicit
formula for \mbox{$\omega_{\diffeo}(\diffeo ')\op{\circ}\sigma_{i}\in
C^{\infty}(\cl{V}_{i},K)$} to show the smoothness of \mbox{$\omega_{\diffeo}$}. Let \mbox{$O_{\diffeo }\se O$} be an open identity
neighbourhood with \mbox{$\diffeo \op{\circ}O_{\diffeo}\op{\circ
}\diffeo^{-1}\se O$} and denote \mbox{$\ol{\diffeo '}=\diffeo \op{\circ }
\diffeo '\op{\circ } \diffeo ^{-1}$} for \mbox{$\diffeo' \in O_{\diffeo
}$}. Fix \mbox{$\diffeo '$} and \mbox{$x\in\cl{V}_{i}$}. Proceeding as in Remark
\ref{rem:valuesOfTheLiftedDiffeomorphismInTermsOfLocalData} we find
\mbox{$j_{\ell},\dots ,j_{1}$} such that
\[
S(\ol{\diffeo'})^{-1}(\sigma_{i}(x))=\sigma_{i_{\ell}}(\ol{\diffeo'
}^{-1}(x))\cdot
k_{j_{\ell}j_{\ell-1}}((\ol{\diffeo}_{j_{\ell-1}})^{-1}\op{\circ}\ldots
\op{\circ}\left(\ol{\diffeo}_{j_{1}}\right)^{-1}(x))\cdot \ldots\cdot
k_{j_{1}i}(x).
\]
Furthermore, let \mbox{$j'_{1}$} be minimal such that 
\[
\big(F_{M}^{-1}\op{\circ}S(\ol{\diffeo '})^{-1}_{M}\big)(x)=\diffeo
^{-1}\op{\circ}\diffeo'^{-1}(x)\in V_{j'_{1}} 
\]
and let \mbox{$U_{x}$} be an open neighbourhood of \mbox{$x$} (w.l.o.g. such that
\mbox{$\cl{U}_{x}\se \cl{V}_{i}$} is a manifold with corners) such that
\mbox{$\ol{\diffeo'}^{-1}(\cl{U}_{x})\se V_{j_{\ell}}$} and
\mbox{$\diffeo^{-1}\op{\circ}\diffeo'^{-1}(\cl{U}_{x})\se V_{j'_{1}}$}. Since
\mbox{$F_{M}=\diffeo$} and
\[
F^{-1}(\sigma_{j_{\ell}}(\ol{\diffeo '}^{-1}(x')))\in
\sigma_{j'_{1}}({\diffeo}^{-1} \op{\circ}\diffeo'^{-1}(x'))\text{ for }
x'\in U_{x}
\]
we have 
\[
F^{-1}\big(\sigma_{j_{\ell}}(\ol{\diffeo '}^{-1}(x'))\big)=
\sigma_{j'_{1}}(\diffeo^{-1}\op{\circ }\diffeo'^{-1}(x'))\cdot
k_{F,x,\diffeo '} ( x') \text{ for } x'\in U_{x},
\]
for some smooth function \mbox{$k_{F,x,\diffeo '}:U_{x}\to K$}. In fact, we
have 
\[
k_{F,x,\diffeo
'}(x)=k_{\sigma_{j'_{1}}}(F^{-1}(\sigma_{j_{\ell}}(\ol{\diffeo
'}^{-1}(x)))).
\]
After possibly shrinking \mbox{$U_{x}$}, a construction as in Remark
\ref{rem:choiceOfLocalTrivialisations} shows that
\mbox{$\left.k_{\sigma_{j'_{1}}}\op{\circ}F^{-1}\op{\circ}
\sigma_{j_{\ell}}\right|_{\cl{U}_{x}}$} extends to a smooth
function on \mbox{$M$}. Thus
\mbox{$\left.k_{F,x,\diffeo' }\right|_{\cl{U}_{x}}\in
C^{\infty}(\cl{U}_{x},K)$} depends smoothly on \mbox{$\diffeo '$} for fixed
\mbox{$x$}.

Accordingly, we find \mbox{$j'_{2},\dots ,j'_{\ell'}$} and a smooth function
\mbox{$k'_{F,x,\diffeo '}:\cl{U}_{x}\to K$} (possibly after shrinking
\mbox{$U_{x}$}), depending smoothly on \mbox{$\diffeo$} such that
\begin{align}\label{eqn:valuesOfTheLiftedDiffeomorphismInSmoothCocycle}
 \omega_{\diffeo}(\diffeo ')(\sigma_{i}(x))=
 \sigma_{i}(x)\cdot \big[ k'_{F,x,\diffeo '}(x)
&\cdot 
 k_{j'_{\ell'}j'_{\ell'-1}}(\diffeo (x))
 \cdot \ldots\cdot 
 k_{j'_{2}j'_{1}}(\diffeo '^{-1}\op{\circ}\diffeo^{-1}(x))\cdot 
 k_{F,x,\diffeo '}(x)\notag\\
&\cdot k_{j_{\ell}j_{\ell-1}}(\diffeo '(x))
 \cdot \ldots\cdot 
 k_{j_{1}i}(x) \big].
\end{align}

Denote the element in brackets on the right hand side by
\mbox{$\kappa_{x,\diffeo '}$}. Since we will not vary \mbox{$F$} and \mbox{$\diffeo$} in
the sequel, we suppressed the dependence of \mbox{$\kappa_{x,\diffeo '}$} on
them. By continuity (cf.\ Remark
\ref{rem:valuesOfTheLiftedDiffeomorphismInTermsOfLocalData}), we find
open neighbourhoods \mbox{$O_{\diffeo' }$} and \mbox{$U'_{x}$} of \mbox{$\diffeo '$} and
\mbox{$x$} (w.l.o.g. such that \mbox{$\cl{U'}_{x}\se\cl{V}_{i}$} is a manifold with
corners) such that
\eqref{eqn:valuesOfTheLiftedDiffeomorphismInSmoothCocycle} defines
\mbox{$\omega_{\diffeo}(h')(\sigma_{i}(x'))$} for all \mbox{$h'\in O_{\diffeo '}$}
and \mbox{$x'\in \cl{U}_{x}$}.  Then \mbox{$\kappa_{x_{1},\diffeo
'}=\kappa_{x_{2},\diffeo '}$} on \mbox{$\cl{U}_{x_{1}}\cap \cl{U}_{x_{2}}$},
finitely many \mbox{$U_{x_{1}},\dots ,U_{x_{m}}$} cover \mbox{$\cl{V}_{i}$} and
since the gluing and restriction maps from Lemma
\ref{lem:restrictionMapForCurrentGroupIsSmooth} and Proposition
\ref{prop:gluingLemmaForCurrentGroup} are smooth,
\[
\omega_{\diffeo}(\diffeo ')\op{\circ }\sigma_{i}=
\glue(
\left.\kappa_{x_{1},\diffeo '}\right|_{\cl{U}_{x_{1}}}
,\dots ,
\left.\kappa_{x_{m},\diffeo '}\right|_{\cl{U}_{x_{m}}}
)
\]
shows that \mbox{$\omega_{\diffeo}(\diffeo ')\op{\circ }\sigma_{i}$} depends
smoothly on \mbox{$\diffeo '$}.
\end{prf}

Before coming to the main result of this section we give a description
of the image of \mbox{$\Diff(M)_{\cP}:=Q(\Aut(\cP))$} in terms of \mbox{$\cP$},
without referring to \mbox{$\Aut(\cP)$}.

\begin{remark}\label{rem:alternativeDescriptionOfDIFF_P} Let \mbox{$Q\from \Aut(\cP)\to \Diff(M)$}, \mbox{$F\mapsto F_{M}$} be the
homomorphism from Definition
\ref{def:bundleAutomorphismsAndGaugeTransformations}. If \mbox{$\diffeo \in
\Diff(M)_{\cP}$}, then there exists an \mbox{$F\in \Aut(\cP)$} that covers
\mbox{$\diffeo$}. Hence the commutative diagram
\[
\begin{CD}
\diffeo^{*}(P)@>\diffeo_{\cP}>> P @>F^{-1}>> P\\
@V\diffeo^{*}(\pi)VV @V\pi VV @V\pi VV \\
M @>\diffeo >> M @>\diffeo^{-1}>> M
\end{CD}
\]
shows that \mbox{$\diffeo^{*}(\cP)$} is equivalent to \mbox{$\cP$}. On the other hand, if
\mbox{$\cP\sim\diffeo^{*}(\cP)$}, then the commutative diagram
\[
\begin{CD}
P@>\sim >> \diffeo^{*}(P)@>\diffeo_{\cP}>> P\\
@V\pi VV @V\diffeo^{*}(\pi)VV @V\pi VV \\
M@= M @>\diffeo >> M 
\end{CD}
\]
shows that there is an \mbox{$F\in \Aut(\cP)$} covering \mbox{$\diffeo$}. Thus
\mbox{$\Diff (M)_{\cP}$} consists of those diffeomorphisms preserving the
equivalence class of \mbox{$\cP$} under pull-backs. This shows also that
\mbox{$\Diff (M)_{\cP}$} is open because homotopic maps yield equivalent
bundles. It thus is contained in \mbox{$\Diff(M)_{0}$}.

Note, that it is not possible to say what \mbox{$\Diff (M)_{\cP}$} is in
general, even in the case of bundles over \mbox{$M=\bS^{1}$}. In fact, we
then have \mbox{$\pi_{0}(\Diff (\bS^{1}))\cong \Z_{2}$} (cf.\
\cite{milnor84}), and the component of \mbox{$\Diff(\bS^{1})$}, which does
not contain the identity, are precisely the orientation reversing
diffeomorphisms on \mbox{$\bS^{1}$}. It follows from the description of
equivalence classes of principal bundles over \mbox{$\bS^{1}$} by
\mbox{$\pi_{0}(K)$} that pulling back the bundle along a orientation
reversing diffeomorphism inverts a representing element for the bundle
in \mbox{$K$}. Thus we have \mbox{$\diffeo^{*}(\cP_{k})\cong \cP_{k^{-1}}$} for
\mbox{$\diffeo \notin \Diff (\bS^{1})_{0}$}. If \mbox{$\pi_{0}(K)\cong \Z_{2}$},
then \mbox{$\cP_{k^{-1}}$} and \mbox{$\cP_{k}$} are equivalent because
\mbox{$[k]=[k^{-1}]$} in \mbox{$\pi_{0}(K)$} and thus \mbox{$\diffeo \in \Diff
(\bS^{1})_{\cP_{k}}$} and \mbox{$\Diff (\bS^{1})_{\cP_{k}}=\Diff
(\bS^{1})$}. If \mbox{$\pi_{0}(K)\cong \Z_{3}$}, then \mbox{$\cP_{k}$} and
\mbox{$\cP_{k^{-1}}$} are \emph{not} equivalent because \mbox{$[k]\neq [k^{-1}]$} in
\mbox{$\pi_{0}(K)$} and thus \mbox{$\diffeo \notin \Diff (\bS^{1})_{\cP_{k}}$} and
\mbox{$\Diff (\bS^{1})_{\cP_{k}}=\Diff (\bS^{1})_{0}$}.
\end{remark}
\begin{theorem}[\mbox{$\boldsymbol{\Aut(\cP)}$} as an extension of
\mbox{$\boldsymbol{\Diff(M)_{\cP}}$} by \mbox{$\boldsymbol{\Gau(\cP)}$}]
\label{thm:automorphismGroupIsLieGroup} Let \mbox{$\cP$} be a smooth
principal \mbox{$K$}-bundle over the closed compact manifold \mbox{$M$}. If \mbox{$\cP$}
has the property SUB, then \mbox{$\Aut(\cP)$} carries a Lie group structure
such that we have an extension of smooth Lie groups
\begin{align}\label{eqn:extensionOfGauByDiff2} 
\Gau(\cP)\hookrightarrow \Aut(\cP)\xrightarrow{Q}
\Diff(M)_{\cP},
\end{align}
where \mbox{$Q:\Aut(\cP)\to \Diff(M)$} is the homomorphism from Definition
\ref{def:bundleAutomorphismsAndGaugeTransformations} and
\mbox{$\Diff(M)_{\cP}$} is the open subgroup of \mbox{$\Diff(M)$} preserving
the equivalence class of \mbox{$\cP$} under pull-backs.
\end{theorem}

\begin{prf}

We identify \mbox{$\Gau(\cP)$} with \mbox{$C^{\infty}(P,K)^{K}$} and extend \mbox{$S$} to a
(possibly non-continuous) section \mbox{$S\from\Diff(M)_{\cP}\to \Aut(\cP)$}
of \mbox{$Q$}.  Now the preceding lemmas show that \mbox{$(T,\omega)$} is a smooth
factor system \cite[Proposition II.8]{neeb06nonAbelianExtensions},
which yields the assertion.
\end{prf}

\begin{proposition}\label{prop:actionOfAutomorphismGroupOnBundleIsSmooth}
In the setting of the previous theorem, the natural action
\[
\Aut(\cP)\times P\to P,\quad (F,p)\mapsto F(p)
\]
is smooth.
\end{proposition}

\begin{prf}
First we note the \mbox{$\Gau(\cP)\cong C^{\infty}(P,K)^{K}$} acts smoothly
on \mbox{$P$} by \mbox{$(\gamma ,p)\mapsto p\cdot \gamma (p)$}. Let \mbox{$O\se \Diff(M)$}
be the neighbourhood from Remark
\ref{rem:chartsForDiffeomorphismGroups} and \mbox{$S\from O\to \Aut(P)$},
\mbox{$g\mapsto \wt{g_{n}}\op{\circ}\dots \op{\circ}\wt{g_{1\,}}$} be the map
from Definition
\ref{def:sectionFromDiffeomorphismsToBundleAutomorphisms}. Then
\mbox{$\Gau(\cP)\op{\circ} S(O)$} is an open neighbourhood in \mbox{$\Aut(\cP)$} and it
suffices to show that the restriction of the action to this
neighbourhood is smooth. Since
\mbox{$\Gau(\cP)$} acts smoothly on \mbox{$\cP$}, this in turn follows from the
smoothness of the map
\[
R\from O\times P\to P,\quad (\diffeo ,p)\mapsto 
S(\diffeo )(p)=\wt{\diffeo_{n}}\op{\circ}\dots \op{\circ}\wt{\diffeo_{1\,}}(p).
\]
To check the smoothness of \mbox{$R$} it suffices to check that
\mbox{$r_{i}\from O\times P\to P\times O$}, \mbox{$(\diffeo ,p)\mapsto
(\wt{\diffeo_{i}}(p),\diffeo)$} is smooth, because then
\mbox{$R=\pr_{1}\op{\circ}r_{n}\op{\circ}\dots \op{\circ}r_{1}$} is
smooth. Now the explicit formula 
\[
\wt{\diffeo_{i}}(\pi (p))=\left\{\begin{array}{ll}
\sigma_{i}(g_{i}(\pi (p)))\cdot k_{i}(p)
&	\text{ if } p\in \pi^{-1}(U_{i})\\
p & \text{ if } p\in \pi^{-1}(\cl{V}_{i})^{c}
\end{array}\right.
\]
shows that \mbox{$r_{i}$} is smooth on \mbox{$\big(O\times \pi^{-1}(U_{i})\big)\cup
\big(O\times \pi^{-1}(\cl{V_{i}})^{c}\big)=O\times P$}.
\end{prf}

\begin{proposition}\label{lem:automorphismGroupActingOnConnections}
If \mbox{\mbox{$\cP$}} is a finite-dimensional smooth principal
\mbox{\mbox{$K$}}-bundle over the closed compact manifold \mbox{\mbox{$M$}}, then
the action
\[
\Aut(\cP)\times \Omega^{1}(P,\fk) \to \Omega^{1}(P,\fk),\quad
F\mapsto (F^{-1})^{*}A,
\]
is smooth. Since this action preserves the closed subspace \mbox{$\Conn(\cP)$} of
connection \mbox{$1$}-forms of \mbox{$\Omega^{1}(P,\fk)$}, the restricted action
\[
\Aut(\cP)\times \Conn(\cP) \to \Conn(\cP),\quad
F\mapsto (F^{-1})^{*}A
\]
is also smooth.
\end{proposition}

\begin{prf}
As in Proposition \ref{prop:actionOfAutomorphismGroupOnBundleIsSmooth}
it can bee seen that the canonical action \mbox{\mbox{$\Aut(P)\times TP\to
TP$}}, \mbox{\mbox{$F.X_{p}=TF(X_{p})$}} is smooth. Since \mbox{\mbox{$P$}} is
assumed to be finite-dimensional and the topology on
\mbox{\mbox{$\Omega^{1}(P,\fk)$}} is the induced topology from
\mbox{\mbox{$C^{\infty}(TP,\fk)$}}, the assertion now follows from
\cite[Proposition 6.4]{gloeckner02patched}.
\end{prf}

\begin{remark}
Of course, the Lie group structure on \mbox{$\Aut(\cP)$} from Theorem
\ref{thm:automorphismGroupIsLieGroup} depends on the choice of \mbox{$S$} and
thus on the choice of the chart \mbox{$\varphi :O\to \cV(M)$} from Remark
\ref{rem:chartsForDiffeomorphismGroups}, the choice of the
trivialising system from Remark \ref{rem:choiceOfLocalTrivialisations}
and the choice of the partition of unity chosen in the proof of Lemma
\ref{lem:decompositionOfDiffeomorphisms}. 

However, different choices lead to isomorphic Lie group structures on
\mbox{$\Aut(\cP)$} and, moreover to equivalent extensions.  To show this we
show that \mbox{$\id_{\Aut(\cP)}$} is smooth when choosing two different
trivialising systems \mbox{$\cl{\cV}=(\cl{V}_{i},\sigma_{i})_{\In}$} and
\mbox{$\cl{\cV}'=(\cl{V}'_{j},\tau_{j})_{\Jm}$}.

Denote by \mbox{$S\from O\to \Aut(\cP)$} and \mbox{$S'\from O\to \Aut(\cP)$} the
corresponding sections of \mbox{$Q$}. Since 
\[
\Gau(\cP)\op{\circ}S(O)=Q^{-1}(O)=\Gau(\cP)\op{\circ}S'(O)
\]
is an open unit neighbourhood and \mbox{$\id_{\Aut(\cP)}$} is an isomorphism
of abstract groups, it suffices to show that the restriction of
\mbox{$\id_{\Aut(\cP)}$} to \mbox{$Q^{-1}(O)$} is smooth.
Now the smooth structure on \mbox{$Q^{-1}(O)$} induced from \mbox{$S$} and \mbox{$S'$}
is given by requiring
\begin{align*}
Q^{-1}(O)\ni 
&F\mapsto (F\op{\circ }S(F_{M})^{-1},F_{M})\in\Gau(\cP)\times \Diff(M)\\
Q^{-1}(O)\ni 
&F\mapsto (F\op{\circ }S'(F_{M})^{-1},F_{M})\in\Gau(\cP)\times \Diff(M)
\end{align*}
to be diffeomorphisms and we thus have to show that 
\[
O\ni \diffeo
\mapsto S(\diffeo)\op{\circ}S'(\diffeo)^{-1}\in \Gau(\cP)
\]
is smooth. By deriving explicit formulae for
\mbox{$S(\diffeo)\op{\circ}S'(\diffeo)^{-1}(\sigma_{i}(x))$} on a
neighbourhood \mbox{$U_{x}$} of \mbox{$x\in \cl{V}_{i}$}, and \mbox{$O_{\diffeo}$} of
\mbox{$\diffeo \in O$} this follows exactly as in Lemma
\ref{lem:smoothCocycleOnAutomorhismGroup}.
\end{remark}

We have not mentioned regularity so far since it is not
needed to obtain the preceding results. However, it is an important
concept and we shall elaborate on it now.

\begin{proposition}
Let \mbox{$\cP$} be a smooth principal \mbox{$K$}-bundle over the compact manifold
with corners \mbox{$M$}. If \mbox{$K$} is regular, then so is \mbox{$\Gau(\cP)$} and,
furthermore, if \mbox{$M$} is closed, then \mbox{$\Aut(\cP)$} is also regular.
\end{proposition}

\begin{prf}
The second assertion follows from the first, because extensions of
regular Lie groups by regular ones are themselves regular
\cite[Theorem
5.4]{omoriMaedaYoshiokaKobayashi83OnRegularFrechetLieGroups} (cf.\
\cite[Theorem V.1.8]{neeb06}).

Let \mbox{$\cl{\cV}=(\ol{V}_{i},\ol{\sigma_{i}})$} be a smooth closed
trivialising system such that \mbox{$\cP$} has the property SUB with respect
to \mbox{$\cl{\cV}$}. We shall use the regularity of
\mbox{$C^{\infty}(\ol{V}_{i},K)$} to obtain the regularity of \mbox{$\Gau(\cP)$}. If
\mbox{$\xi \from [0,1]\to \gau(\cP)$} is smooth, then this determines smooth
maps \mbox{$\xi_{i}\from [0,1]\to C^{\infty}(\ol{V}_{i},\fk)$}, satisfying
\mbox{$\xi_{i}(t)(m)=\Ad(k_{ij}(m)).\xi_{j}(t)(m)$}. By regularity of
\mbox{$C^{\infty}(\ol{V}_{i},K)$} this determines smooth maps
\mbox{$\gamma_{\xi,i}\from [0,1]\to C^{\infty}(\ol{V}_{i},K)$}.

By uniqueness of solutions of differential equations we see that the
mappings \mbox{\mbox{$t\mapsto \gamma_{\xi ,i}(t)$}} and \mbox{\mbox{$t\mapsto
k_{ij}\cdot \gamma_{\xi ,j}(t)\cdot k_{ji}$}} have to coincide,
ensuring \mbox{\mbox{$\gamma_{\xi ,i}(t)(m)=k_{ij}(m)\cdot
\gamma_{\xi ,j}(t)(m)\cdot k_{ji}(m)$}} for all \mbox{$t\in [0,1]$} and \mbox{$m\in
\ol{V}_{i}\cap \ol{V}_{j}$}.  Thus \mbox{$[0,1]\ni t\mapsto (\gamma_{\xi
,i}(t))_{\In}\in G_{\cl{\cV}}(\cP )$} is a solution of the
corresponding initial value problem and the desired properties follows
from the regularity of \mbox{$C^{\infty}(\ol{V}_{i},K)$}.
\end{prf}

\begin{remark}
A Lie group structure on \mbox{$\Aut(\cP)$} has been considered in
\cite{michorAbbatiCirelliMania89automorphismGroup} in the convenient
setting, and the interest in \mbox{$\Aut(\cP)$} as a symmetry group coupling
the gauge symmetry of Yang-Mills theories and the
\mbox{$\Diff(M)$}-invariance of general relativity is emphasised. Moreover,
it is also shown that \mbox{$\Gau(\cP)$} is a split Lie subgroup of
\mbox{$\Aut(\cP)$}, that
\[
\Gau(\cP)\hookrightarrow \Aut(\cP)\twoheadrightarrow \Diff(M)_{\cP}
\]
is an exact sequence of Lie groups and that the action
\mbox{$\Aut(\cP)\times P\to P$} is smooth. However, the Lie group structure
is constructed out of quite general arguments allowing to give the
space \mbox{$\Hom(\cP,\cP)$} of bundle morphisms a smooth structure and then
to consider \mbox{$\Aut(\cP)$} as an open subset of \mbox{$\Hom(\cP,\cP)$}.

The approach taken in this section is somehow different,
since the Lie group structure on \mbox{$\Aut(\cP)$} is constructed by foot
and the construction provides explicit charts given by charts of
\mbox{$\Gau(\cP)$} and \mbox{$\Diff(M)$}.
\end{remark}

\begin{remark}
The approach to the Lie group structure in this section used detailed
knowledge on the chart \mbox{$\varphi \from O\to \cV(M)$} of the Lie group
\mbox{$\Diff(M)$} from Remark \ref{rem:chartsForDiffeomorphismGroups}. We
used this when decomposing a diffeomorphism into a product of
diffeomorphisms with support in some trivialising subset of \mbox{$M$}. The
fact that we needed was that for a diffeomorphism \mbox{$\diffeo \in O$} we
have \mbox{$\diffeo (m)=m$} if the vector field \mbox{$\varphi (\diffeo )$} vanishes
in \mbox{$m$}. This should also be true for the charts on \mbox{$\Diff (M)$} for
compact manifolds with corners and thus the procedure of this section
should carry over to bundles over manifolds with corners.
\end{remark}

\begin{example}[\mbox{\mbox{$\boldsymbol{\Aut(C^{\infty}_{k}(\bS^{1},\fk))}$}}]
\label{exmp:automorphismGroupOfTwistedLoopAlgebra}Let \mbox{\mbox{$K$}} be a
simple finite-dimensional Lie group, \mbox{\mbox{$K_{0}$}} be compact and
simply connected and \mbox{\mbox{$\cP_{k}$}} be a smooth principal
\mbox{\mbox{$K$}}-bundle over \mbox{\mbox{$\bS^{1}$}}, uniquely determined up to
equivalence by \mbox{$[k]\in \pi_{0}(K)$}. Identifying the twisted loop algebra
\[
C^{\infty }_{k}(\bS^{1},\fk):=\{\eta \in C^{\infty}(\R,\fk):\eta
(x+n)=\Ad(k)^{-n}.\eta (x)\fa x\in\R,n\in\Z\}.
\]
with the gauge algebra of the flat principal bundle \mbox{$\cP_{k}$}, we get
a smooth action of \mbox{\mbox{$\Aut(\cP_{k})$}} on
\mbox{\mbox{$C^{\infty}_{k}(\bS^{1},\fk)$}}, which can also be lifted to the
twisted loop group \mbox{$C^{\infty}_{k}(\bS^{1},K)$}, the affine Kac--Moody
algebra \mbox{$\wh{C^{\infty}_{k}(\bS^{1},\fk)}$} and to the affine
Kac--Moody group \mbox{$\wh{C^{\infty}_{k}(\bS^{1},K)}$} \cite{diss}. Various
results (cf.\ \cite[Theorem 16]{lecomte80}) assert that each
automorphism of \mbox{\mbox{$C^{\infty}_{k}(\bS^{1},\fk)$}} arises in this
way and we thus have a geometric description of
\mbox{\mbox{$\Aut(C^{\infty}_{k}(\bS^{1},\fk))\cong
\Aut(C^{\infty}_{k}(\bS^{1},K)_{0})\cong \Aut(\cP_{k})$}} for
\mbox{$C^{\infty}_{k}(\bS^{1},K)_{0}$} is simply connected. Furthermore,
this also leads to topological information on
\mbox{\mbox{$\Aut(C^{\infty}_{k}(\bS^{1},\fk))$}}, since we get a long exact
homotopy sequence
\begin{multline}\label{eqn:longExactHomotopySequenceForAutomorphismGroup}
\dots 
\to 
\pi_{n+1}(\Diff(\bS^{1}))
\xrightarrow{\delta_{n+1}}
\pi_{n}(C^{\infty}_{k}(\bS^{1},K))
\to 
\pi_{n}(\Aut(\cP_{k}))\\
\to
\pi_{n}(\Diff(\bS^{1}))
\xrightarrow{\delta_{n}}
\pi_{n-1}(C^{\infty}_{k}(\bS^{1},K))
\to
\dots 
\end{multline}
induced by the locally trivial bundle
\mbox{\mbox{$\Gau(\cP_{k})\hookrightarrow
\Aut(\cP_{k})\xrightarrow{q} \Diff (\bS^{1})_{\cP_{k}}$}} and
the isomorphisms \mbox{\mbox{$\Gau(\cP_{k})\cong C^{\infty}_{k}(\bS^{1},K)$}}
and \mbox{\mbox{$\Aut(\cP_{k})\cong \Aut(C^{\infty}_{k}(\bS^{1},\fk))$}}.
E.g., in combination with
\begin{align}\label{eqn:homotopyGroupsForDiffeomorphismGroup}
\pi_{n}(\Diff(\bS^{1}))\cong\left\{\begin{array}{ll}
\Z_{2} &	\text{ if }n=0\\
\Z &	\text{ if }n=1\\
0 &	\text{ if }n\geq 2
\end{array} \right.
\end{align}
(cf.\ \cite{milnor84}), one obtains information on
\mbox{\mbox{$\pi_{n}(\Aut(\cP_{k}))$}}. In fact, consider the exact sequence
\begin{multline*}
0
\to
\pi_{1}(C^{\infty}_{k}(\bS^{1},K))
\to
\pi_{1}(\Aut(\cP_{k}))
\to\underbrace{\pi_{1}(\Diff (M))}_{\cong \Z}
\xrightarrow{\delta_{1}}
\pi_{0}(C^{\infty}_{k}(\bS^{1},K))\\
\to \pi_{0}(\Aut(\cP_{k}))\xrightarrow{\pi_{0}(q)} 
\pi_{0}(\Diff(\bS^{1})_{\cP_{k}}) 
\end{multline*}
induced by \eqref{eqn:longExactHomotopySequenceForAutomorphismGroup}
and \eqref{eqn:homotopyGroupsForDiffeomorphismGroup}.
Since \mbox{\mbox{$\pi_{1}(C^{\infty}_{k}(\bS^{1},K))$}} vanishes, this
implies \mbox{\mbox{$\pi_{1}(\Aut(\cP_{k}))\cong \Z$}}. A generator of
\mbox{$\pi_{1}(\Diff (\bS^{1}))$} is \mbox{$\id_{\bS^{1}}$}, which lifts to a
generator of \mbox{$\pi_{1}(\Aut(\cP_{k}))$}. Thus the connecting
homomorphism \mbox{$\delta_{1}$} vanishes.  The argument from Remark
\ref{rem:alternativeDescriptionOfDIFF_P} shows precisely that
\mbox{\mbox{$\pi_{0}(\Diff (\bS^{1})_{\cP_{k}})\cong \Z_{2}$}} if and only if
\mbox{\mbox{$k^{2}\in K_{0}$}} and that \mbox{$\pi_{0}(q)$} is surjective.  We thus
end up with an exact sequence
\[
\operatorname{Fix}_{\pi_{0}(K)} ([k])\to
\pi_{0}(\Aut(\cP_{k}))\twoheadrightarrow \left\{\begin{array}{ll}
\Z_{2} & \text{ if }k^{2}\in K_{0}\\
\1 & \text{ else.}
\end{array} \right.
\]
Since
\eqref{eqn:homotopyGroupsForDiffeomorphismGroup} implies that
\mbox{\mbox{$\Diff(\bS^{1})_{0}$}} is a \mbox{\mbox{$K(1,\Z)$}}, we also have
\mbox{\mbox{$\pi_{n}(\Aut(\cP_{k}))\cong
\pi_{n}(C^{\infty}_{k}(\bS^{1},K))$}} for \mbox{\mbox{$n\geq 2$}}.
\end{example}

\begin{remark}
The description of \mbox{\mbox{$\Aut(C^{\infty}_{k}(\bS^{1},\fk))$}} in Example
\ref{exmp:automorphismGroupOfTwistedLoopAlgebra} should arise out of a
general principle, describing the automorphism group for gauge
algebras of flat bundles, i.e., of bundles of the
form
\[
\cP_{\varphi}=\wt{M}\times K/\sim \;\text{ where }\;
(m,k)\sim (m\cdot d,\varphi^{-1}(d)\cdot k).
\]
Here \mbox{\mbox{$\varphi \from \pi_{1}(M)\to K$}} is a homomorphism and \mbox{\mbox{$\wt{M}$}} is
the simply connected cover of \mbox{\mbox{$M$}}, on which \mbox{\mbox{$\pi_{1}(M)$}} acts
canonically. Then
\[
\gau(\cP)\cong C^{\infty}_{\varphi}(M,\fk):=\{\eta \in
C^{\infty}(\wt{M},\fk): \eta (m\cdot d)=\Ad (\varphi (d))^{-1}.\eta
(m)\}.
\]
and this description should allow to reconstruct gauge transformations
and diffeomorphisms out of the ideals of
\mbox{\mbox{$C^{\infty}_{\varphi}(M,\fk)$}} (cf.\ \cite{lecomte80}).
\end{remark}

\begin{problem}(cf.\ \cite[Problem IX.5]{neeb06}) Let
\mbox{\mbox{$\cP_{\varphi}$}} be a (flat) principal \mbox{\mbox{$K$}}-bundle over
the closed compact manifold \mbox{\mbox{$M$}}. Determine the automorphism
group \mbox{\mbox{$\Aut(\gau (\cP))$}}. In which cases does it coincide with
\mbox{\mbox{$\Aut(\cP)$}} (the main point here is the surjectivity of the
canonical map \mbox{\mbox{$\Aut(\cP)\to \Aut(\gau(\cP))$}}).
\end{problem}

\begin{remark}
In some special cases, the extension \mbox{$\Gau(\cP)\hookrightarrow
\Aut(\cP)\twoheadrightarrow \Diff(M)_{\cP}$} from Theorem
\ref{thm:automorphismGroupIsLieGroup} splits. This is the case for
trivial bundles and for bundles with abelian structure group \mbox{$K$}, but
also for frame bundles, since we then have a natural homomorphism
\mbox{$\Diff(M)\to \Gau(\cP)$}, \mbox{$\diffeo \mapsto d\diffeo$}. However, it would
be desirable to have a characterisation of the bundles, for which
this extension splits.
\end{remark}

\begin{problem}(cf.\ \cite[Problem V.5]{neeb06})
Find a characterisation of those principal \mbox{$K$}-bundles \mbox{$\cP$} for which
the extension \mbox{$\Gau(\cP)\hookrightarrow \Aut(\cP)\twoheadrightarrow
\Diff(M)_{\cP}$} splits on the group level.
\end{problem}

\appendix

\section{Appendix: Differential calculus on spaces of mappings}

\begin{definition}
\label{def:diffcalcOnLocallyConvexSpaces} Let \mbox{$X$} and \mbox{$Y$} be a locally
convex spaces and \mbox{$U\se X$} be open. Then \mbox{$f\from U\to Y$} is
\emph{differentiable}  or \emph{\mbox{$C^{1}$}} if it is
continuous, for each \mbox{$v\in X$} the differential quotient
\[
df (x).v:=\lim_{h\to 0}\frac{f (x+hv)-f (x)}{h}
\]
exists and if the map \mbox{$df\from U\times X\to Y$} is continuous.  If
\mbox{$n>1$} we inductively define \mbox{$f$} to be \emph{\mbox{$C^{n}$}} if it is \mbox{$C^{1}$}
and \mbox{$df$} is \mbox{$C^{n-1}$} and to be \mbox{$C^{\infty}$} or
\emph{smooth}  if it is \mbox{$C^{n}$}. We
say that \mbox{$f$} is \mbox{$C^{\infty}$} or \emph{smooth} if \mbox{$f$} is \mbox{$C^{n}$} for
all \mbox{$n\in\N_{0}$}. We denote the corresponding spaces of maps by
\mbox{$C^{n}(U,Y)$} and \mbox{$C^{\infty}(U,Y)$}. \end{definition}  \begin{definition}\label{def:diffcalcOnNonopenDomains}
Let \mbox{$X$} and \mbox{$Y$} be locally convex spaces, and let \mbox{$U\se X$} be a set
with dense interior. Then \mbox{$f\from U\to Y$} is
\emph{differentiable}  or
\emph{\mbox{$C^{1}$}} if it is continuous, \mbox{$f_{\op{int}}:=\left.f
\right|_{\op{int}(U)}$} is \mbox{$C^{1}$} and the map
\[
d\left(f_{\op{int}}\right):\op{int}(U)\times X \to Y,\;\; (x,v)\mapsto
d\left(f_{\op{int}} \right)(x).v
\]
extends to a continuous map on \mbox{$U\times X$}, which is
called the \emph{differential} \mbox{$df$} of \mbox{$f$}. If
\mbox{$n>1$} we inductively define \mbox{$f$} to be \emph{\mbox{$C^{n}$}} if it is \mbox{$C^{1}$}
and \mbox{$df$} is \mbox{$C^{n-1}$}. We say that \mbox{$f$} is \mbox{$C^{\infty}$} or
\emph{smooth}  if \mbox{$f$} is
\mbox{$C^{n}$} for all \mbox{$n\in \N_{0}$}. We denote the corresponding spaces of
maps by \mbox{$C^{n}(U,Y)$} and
\mbox{$C^{\infty}(U,Y)$}.
\end{definition}

\begin{definition}\label{def:lieGroup}
From the definition above, the notion of a \textit{Lie group} is
clear. It is a group which is a smooth manifold modelled on a locally
convex space such that the group operations are smooth.  Moreover, the
notion of a finite-dimensional manifold with corners is clear, i.e., a
smooth (in the sense of Definition \ref{def:diffcalcOnNonopenDomains})
manifold modelled on $\R^{n}_{+}:=\{(x_{1},\dots ,x_{n})\in
\R^{n}:x_{i}\geq 0\}$ (cf.\ \cite{michor80} and \cite{smoothExt}).
\end{definition}

\begin{proposition}
\label{prop:localDescriptionsOfLieGroups} Let \mbox{\mbox{$G$}} be a group with a
locally convex manifold structure on some subset \mbox{\mbox{$U\se G$}} with \mbox{\mbox{$e\in
U$}}. Furthermore, assume that there exists \mbox{\mbox{$V\se U$}} open such that
\mbox{\mbox{$e\in V$}}, \mbox{\mbox{$VV\se U$}}, \mbox{\mbox{$V=V^{-1}$}} and
\begin{itemize}
\item [i)] \mbox{\mbox{$V\times V\to U$}}, \mbox{\mbox{$(g,h)\mapsto gh$}} is smooth,
\item [ii)] \mbox{\mbox{$V\to V$}}, \mbox{\mbox{$g\mapsto g^{-1}$}} is smooth,
\item [iii)] for all \mbox{\mbox{$g\in G$}}, there exists an open unit neighbourhood
\mbox{\mbox{$W\se U$}} such that \mbox{\mbox{$g^{-1}Wg\se U$}} and the map \mbox{\mbox{$W\to U$}}, \mbox{\mbox{$h\mapsto
g^{-1}hg$}} is smooth.
\end{itemize}
Then there exists a unique locally convex manifold structure on \mbox{\mbox{$G$}}
which turns \mbox{\mbox{$G$}} into a Lie group, such that \mbox{\mbox{$V$}} is an open submanifold
of \mbox{\mbox{$G$}}.
\end{proposition}

\begin{definition}\label{def:exponentialfunction}
Let \mbox{\mbox{$G$}} be a locally convex Lie group. The group \mbox{\mbox{$G$}} is said to have
an \emph{exponential function} if for each \mbox{\mbox{$x \in \fg$}} the initial
value problem
\[
\gamma (0)=e,\quad \gamma '(t)=T\lambda_{\gamma (t)}(e).x
\]
has a solution \mbox{\mbox{$\gamma_{x}\in C^{\infty} (\R,G)$}} and the
function
\[
\exp_{G}:\fg\to G,\;\;x \mapsto \gamma_x (1)
\]
is smooth. Furthermore, if there exists a zero neighbourhood \mbox{\mbox{$W\se
\fg$}} such that \mbox{\mbox{$\left.\exp_{G}\right|_{W}$}} is a diffeomorphism onto
some open unit neighbourhood of \mbox{\mbox{$G$}}, then \mbox{\mbox{$G$}} is said to be
\emph{locally exponential}.
\end{definition}

\begin{lemma}
\label{lem:interchangeOfActionsOnGroupAndAlgebra} If \mbox{\mbox{$G$}} and \mbox{\mbox{$G'$}} are
locally convex Lie groups with exponential function, then for each
morphism \mbox{\mbox{$\alpha :G\to G'$}} of Lie groups and the induced morphism
\mbox{\mbox{$d\alpha (e):\fg \to \fg'$}} of Lie algebras, the diagram
\[
\begin{CD}
G @>\alpha >>G'\\
@AA\exp_{G}A  @AA\exp_{G'}A\\
\fg@>d\alpha (e)>> \fg'
\end{CD}
\]
commutes.
\end{lemma}

\begin{remark}
\label{rem:banachLieGroupsAreLocallyExponential} The Fundamental
Theorem of Calculus for locally convex spaces (cf.\ \cite[Theorem
1.5]{gloeckner02b}) yields that a locally convex Lie group \mbox{\mbox{$G$}} can
have at most one exponential function (cf.\ \cite[Lemma
II.3.5]{neeb06}). If \mbox{\mbox{$G$}} is a Banach-Lie group (i.e., \mbox{\mbox{$\fg$}} is a
Banach space), then \mbox{\mbox{$G$}} is locally exponential due to the
existence of solutions of differential equations, their smooth
dependence on initial values \cite[Chapter IV]{lang99} and the Inverse
Mapping Theorem for Banach spaces \cite[Theorem I.5.2]{lang99}. In
particular, each finite-dimensional Lie group is locally exponential.
\end{remark}

\begin{definition}\label{def:topologiesOnSpacesOfMappings} If \mbox{$X$} is a Hausdorff space
and \mbox{$Y$} is a topological spaces, then the \emph{compact-open
topology}  on
the space of continuous functions is defined as the topology generated
by the sets of the form
\[
\lfloor C,W\rfloor:=\{f\in C(X,Y): f(C)\se W\}, 
\] where \mbox{$C$} runs over all compact subsets of \mbox{$X$} and \mbox{$W$} runs over all
open subsets of \mbox{$Y$}. We write
\mbox{$C(X,Y)_{c}$} for the space \mbox{$C(X,Y)$} endowed with the compact-open topology.

If \mbox{$G$} is a topological group, then \mbox{$C(X,G)$} is a group with respect
to pointwise group operation. Furthermore, the topology of compact
convergence coincides with the compact-open topology \cite[Theorem
X.3.4.2]{bourbakiTop} and thus \mbox{$C(X,G)_{c}$} is again a topological
group.  A basis of unit neighbourhoods of this topology is given by
\mbox{$\lfloor C,W\rfloor$}, where \mbox{$C$} runs over all compact subsets of \mbox{$X$}
and \mbox{$W$} runs over all open unit neighbourhoods of \mbox{$G$}.  If \mbox{$X$} itself
is compact, then this basis is already given by \mbox{$\lfloor X,W\rfloor$},
where \mbox{$W$} runs over all unit neighbourhoods of \mbox{$G$}.

If \mbox{$Y$} is a locally convex space, then \mbox{$C(X,Y)$} is a vector space with
respect to pointwise operations. The preceding discussion implies
that addition is continuous and scalar multiplication is also
continuous. Since its topology is induced by the seminorms
\[
p_{C}:C(X,Y)\to \K,\quad f\mapsto \sup\nolimits_{x\in C}\{p(f(x))\},
\]
where \mbox{$C$} runs over all compact subsets of \mbox{$X$} and \mbox{$p$} runs over all
seminorms, defining the topology on \mbox{$Y$}, we see that \mbox{$C(X,Y)_{c}$} is again
locally convex.

If \mbox{$M$} and \mbox{$N$} are manifolds with corners as in \cite{smoothExt},
then every smooth map \mbox{$f:M\to N$} defines a sequence of continuous map
\mbox{$T^{n}f:T^{n}M\to T^{n}N$} on the iterated tangent bundles. We thus
obtain an inclusion
\[
C^{\infty}(M,N)\hookrightarrow \prod_{n=0}^{\infty}C(T^{n}M,T^{n}M)_{c},
\quad f\mapsto (T^{n}f)_{n\in\N}
\]
and we define the
\emph{\mbox{$C^{\infty}$}-topology}  on \mbox{$C^{\infty}(M,N)$} to be the initial
topology induced from this inclusion. For a locally convex space \mbox{$Y$}
we thus get a locally convex vector topology on \mbox{$C^{\infty}(M,Y)$}.

If \mbox{$\cE =(Y,\xi :E\to X)$} is a continuous vector bundle and
\mbox{$S_{c}(\cE)$} is the set of continuous sections, then we have an
inclusion \mbox{$S_{c}(\cE)\hookrightarrow C(X,E)$} and we thus obtain
a topology on \mbox{$S_{c}(\cE )$}.  If \mbox{$\cE$} is also smooth, then we have an
inclusion \mbox{$S(\cE)\hookrightarrow C^{\infty}(M,E)$}, inducing a
topology \mbox{$S(\cE)$}, which we also call \mbox{$C^{\infty}$}-topology.
\end{definition}

\begin{remark}
\label{rem:alternativeDescriptionOfTopologyOnSpaceOfSmoothMappings} If
\mbox{$M$} is a manifold with corners and \mbox{$Y$} is a locally convex space, then
we can describe the \mbox{$C^{\infty}$}-topology on \mbox{$C^{\infty}(M,Y)$}
alternatively as the initial topology with respect to the inclusion
\[
C^{\infty}(M,Y)\hookrightarrow \prod_{n=0}^{\infty}C(T^{n}M,Y),\quad 
f\mapsto (d^{n}f)_{n\in \N},
\]
where \mbox{$d^{n}f=\pr_{2^{n}}\op{\circ}T^{n}f$}. In fact, we have
\mbox{$Tf=(f,df)$} and we can inductively write \mbox{$T^{n}f$} in terms of
\mbox{$d^{l}f$} for \mbox{$l\leq n$}. This implies for a map into \mbox{$C^{\infty}(M,Y)$}
that its composition with each \mbox{$d^{n}$} is
continuous if and only if its composition with all \mbox{$T^{n}$} is
continuous. Because the initial topology is characterised by this
property, the topologies coincide.
\end{remark}

\begin{definition}
If \mbox{$\cE=(Y,\xi \from E\to M)$} is a smooth vector bundle and \mbox{$p\in
\N_{0}$}, then a \emph{\mbox{$\cE$}-valued \mbox{$p$}-form}  on \mbox{$M$} is a function \mbox{$\form$} which
associates to each \mbox{$m\in M$} a \mbox{$p$}-linear alternating map
\mbox{$\form_{m}\from (T_{m}M)^{p}\to E_{m}$} such that in local coordinates
the map 
\[
(m,X_{1,m},\dots ,X_{p,m})\mapsto \form_{m}(X_{1,m},\dots
,X_{p,m})
\]
is smooth. We denote by
\[
\Omega^{p}(M,\cE ):=\{\form :\bigcup_{m\in M}(T_{m}M)^{p}\to E:
\form \text{ is a }\cE \text{ valued }p\text{-form on }M\}
\]
the space of \mbox{$\cE$}-valued \mbox{$p$}-forms on \mbox{$M$} which has a
canonical vector space structure induced from pointwise operations.
\end{definition}

\begin{remark}
If \mbox{$\cE =(Y,\xi \from E\to M)$} is a smooth vector bundle over the
finite-dimensional manifold \mbox{$M$}, then each \mbox{$\cE$}-valued \mbox{$p$}-form
\mbox{$\form$} maps vector fields \mbox{$X_{1},\dots ,X_{p}$} to a smooth section
\mbox{$\form.(X_{1},\dots ,X_{p}):=\form \circ (X_{1}\times\dots\times
X_{p})$} in \mbox{$S(\cE)$}, which is \mbox{$C^{\infty}(M,\R)$}-linear by
definition.  Conversely, any alternating \mbox{$C^{\infty}(M)$}-linear map
\mbox{$\form \from \cV^{p}(M)\to S(\cE)$} determines uniquely an element of
\mbox{$\Omega^{p}(M,\cE)$} by setting
\[
\form_{m}(X_{1,m},\dots,X_{n,m}):=\form (\wt{X}_{1},\dots
,\wt{X}_{p})(m),
\]
where \mbox{$\wt{X_{i}}$} is an extension of \mbox{$X_{i,m}$} to a smooth vector
field. That \mbox{$\form_{m}(X_{1,m},\dots ,X_{p,m})$} does not depend on the
choice of this extension follows from the \mbox{$C^{\infty}(M,\R)$}-linearity of
\mbox{$\form$}, if one expands different choices in terms of basis vector
fields. Note that the assumption on \mbox{$M$} to be finite-dimensional is
crucial for this argument.
\end{remark}

We now consider the continuity properties of some very basic maps,
i.e., restriction maps and gluing maps. These maps we shall encounter
often in the sequel.

\begin{lemma}\label{lem:restrictionMapIsContinuous} If \mbox{$\cE$} is a smooth vector bundle over \mbox{$M$} and \mbox{$U\se M$} is open and
\mbox{$\cE_{U}=\left.\cE\right|_{U}$} is the restricted vector bundle, then
the restriction map \mbox{$\res_{U}:S(\cE)\to S(\cE_{U})$}, \mbox{$\sigma \mapsto
\left.\sigma\right|_{U}$} is continuous.
If, moreover, \mbox{$\cl{U}$} is a manifold with corners, then the
restriction map \mbox{$\res_{\cl{U}}:S(\cE)\to S(\cE_{\cl{U}})$}, \mbox{$\sigma
\mapsto \left.\sigma\right|_{\cl{U}}$} is continuous.
\end{lemma}

\begin{proposition}\label{prop:spaceOfSectionOfAVectorBund;eInLocallyConvex}
If \mbox{$\cE$} is a smooth vector bundle over the finite-dimensional manifold
with corners \mbox{$M$} and \mbox{$S(\cE)$} is the vector space of smooth sections
with pointwise operations, then the \mbox{$C^{\infty}$}-topology is a locally
convex vector topology on \mbox{$S(\cE)$}. Furthermore, if \mbox{$(U_{i})_{i\in I}$}
is an open cover of \mbox{$M$} such that each \mbox{$\cl{U}_{i}$} is a
manifold with corners and \mbox{$\cE_{i}:=\left.\cE\right|_{\cl{U}_{i}}$}
denotes the restricted bundle, then the \mbox{$C^{\infty}$}-topology on
\mbox{$S(\cE)$} is initial with respect to
\begin{align}\label{eqn:restrictionMapForVectorBundles}
\res \from S(\cE)\to \prod_{i\in I} S(\cE_{i}),\quad 
\sigma \mapsto (\left.\sigma\right|_{\cl{U}_{i}})_{i\in I}.
\end{align}
\end{proposition}

\begin{corollary}\label{cor:restrictionMapIsSmooth}
The restriction maps \mbox{$\res_{U}$} and \mbox{$\res_{\cl{U}}$} from Lemma
\ref{lem:restrictionMapIsContinuous} are smooth.
\end{corollary}

\begin{proposition}\label{prop:gluingForVectorBundles_CLOSED_Version}
If \mbox{$\cE$} is a smooth vector bundle over the finite-dimensional
manifold with corners \mbox{$M$}, \mbox{$\mathfrak{U}=(U_{i})_{i\in I}$} is an open
cover of \mbox{$M$} such that each \mbox{$\cl{U}_{i}$} is
a manifold with corners and
\mbox{$\cE_{i}:=\left.\cE\right|_{\cl{U}_{i}}$} denotes the restricted
bundle, then
\[
S_{\ol{\mathfrak{U}}}(\cE)=\{(\sigma_{i})_{i\in I}
\in \bigoplus_{i\in I}S(\cE_{i}):
\sigma_{i}(x)=\sigma_{j}(x)\fa x\in \cl{U}_{i}\cap \cl{U}_{j}\}
\] is a closed subspace of \mbox{$\bigoplus_{i\in I}S(\cE_{i})$} and the
\emph{gluing map}
\begin{align}\label{eqn:gluingMapForVectorBundles_COLSED_Verison}
\glue \from S_{\ol{\mathfrak{U}}}(\cE)\to S(\cE),\quad
\glue ((\sigma_{i})_{i\in I})(x)=\sigma_{i}(x)\text{ if }x\in \cl{U}_{i}
\end{align} is inverse to the restriction map
\eqref{eqn:restrictionMapForVectorBundles}.
\end{proposition}

\begin{corollary}\label{cor:gluingForVectorBundles_OPEN_Version}
If \mbox{$\cE$} is a smooth vector bundle over the finite-dimensional
manifold with corners \mbox{$M$}, \mbox{$\mathfrak{U}=(U_{i})_{i\in I}$} is an open
cover of \mbox{$M$} and \mbox{$\cE_{i}:=\left.\cE\right|_{U_{i}}$} denotes the
restricted bundle, then
\[
S_{{\mathfrak{U}}}(\cE)=\{(\sigma_{i})_{i\in I}\in \bigoplus_{i\in
I}S(\cE_{i}): \sigma_{i}(x)=\sigma_{j}(x)\fa x\in U_{i}\cap U_{j}\}
\] is a
closed subspace of \mbox{$\bigoplus_{i\in I}S(\cE_{i})$} and the \emph{gluing
map}
\begin{align}\label{eqn:gluingMapForVectorBundles_OPEN_Version}
\glue \from S_{\mathfrak{U}}(\cE)\to S(\cE),\quad
\glue ((\sigma_{i})_{i\in I})(x)=\sigma_{i}(x)\text{ if }x\in U_{i}
\end{align}
is inverse to the restriction map.
\end{corollary}

\begin{lemma}\label{lem:restrictionMapForCurrentGroupIsSmooth}
If \mbox{$M$} is a compact manifold with corners, \mbox{$K$} is a Lie group and
\mbox{$\cl{U}\se M$} is a manifold with corners, then the restriction
\[
\res : C^{\infty}(M,K)\to C^{\infty}(\cl{U},K),\quad
\gamma \mapsto \left.\gamma\right|_{\cl{U}}
\] is a smooth homomorphism of Lie groups.
\end{lemma}

\begin{proposition}\label{prop:gluingLemmaForCurrentGroup}
Let \mbox{$K$} be a Lie group, \mbox{$M$} be a compact manifold with corners with an
open cover \mbox{$\mathfrak{V}=(V_{1},\dots ,V_{n})$}
such that \mbox{$\cl{\mathfrak{V}}=(\cl{V_{1}},\dots ,\cl{V_{n}})$} is
a cover by manifolds with corners. Then
\[
G_{\cl{\mathfrak{V}}}:=\{(\gamma_{1},\dots ,\gamma_{n})\in
\prod_{i=1}^{n}C^{\infty}(\cl{V}_{i},K):\gamma_{i}(x)=\gamma_{j}(x)\fa x\in
\cl{V}_{i}\cap \cl{V}_{j}\}
\]
is a
closed subgroup of \mbox{$\prod_{i=1}^{n}C^{\infty}(\cl{V}_{i},K)$},
which is a Lie group modelled on the closed subspace
\[
\fg_{\mathfrak{V}}:=\{(\eta_{1},\dots ,\eta_{n})\in
\prod_{i=1}^{n}C^{\infty}(\cl{V}_{i},\fk):\eta_{i}(x)=\eta_{j}(x)\fa x\in
\cl{V}_{i}\cap \cl{V}_{j}\}
\]
of \mbox{$\bigoplus_{i=1}^{n}C^{\infty}(\cl{V}_{i},\fk)$} and the \emph{gluing map}
\[
\glue :G_{\cl{\mathfrak{V}}}\to C^{\infty}(M,K),\quad
\glue (\gamma_{1},\dots ,\gamma_{n})= \gamma_{i}(x)\text{ if }x\in \cl{V}_{i}
\] is an isomorphism of Lie groups.
\end{proposition}

\begin{lemma}
\label{lem:productrule} Let \mbox{\mbox{$M$}}
be a smooth locally convex manifold with corners, \mbox{\mbox{$G$}} be a locally
convex Lie group and \mbox{\mbox{$\lambda:G\times Y\to Y$}} be a smooth linear
action on the locally convex space \mbox{\mbox{$Y$}}. If \mbox{\mbox{$h:M\to G$}} and \mbox{\mbox{$f:M\to Y$}}
are smooth, then we have
\begin{align}\label{eqn:productrule}
d\left(\lambda(h).f \right).X_{m}=\lambda (h).(df.X_{m})+ \dot\lambda
\left(\Ad(h).\delta^{l}(h).X_{m}\right).\left(\lambda (h(m)).f(m)\right)
\end{align}
with \mbox{\mbox{$\lambda(h^{-1}).f:M\to E$}}, \mbox{\mbox{$m\mapsto
\lambda\left(h(m)^{-1}\right).f(m)$}}.  If \mbox{\mbox{$\lambda=\Ad$}} is the adjoint
action of \mbox{\mbox{$G$}} on \mbox{\mbox{$\fg$}}, then we have
\begin{align*}
d\left(\Ad(h).f \right).X_{m}= \Ad(h).(df.X_{m})+
\Ad(h).\left[\delta^{l}(h).X_{m},f(m)\right]
\end{align*}
\end{lemma}

\begin{definition}Let \mbox{\mbox{$\cP=\pfb[X]$}} be a
continuous principal \mbox{\mbox{$K$}}-bundle. If \mbox{\mbox{$(U_{i})_{i\in I}$}} is an open
cover of \mbox{\mbox{$X$}} by trivialising neighbourhoods and \mbox{\mbox{$(\sigma_{i}\from
U_{i}\to P)_{i\in I}$}} is a collection of continuous sections, then the
collection \mbox{\mbox{$\mathcal{U}=(U_{i},\sigma_{i})_{i\in I}$}}
is called an \emph{continuous open trivialising system} of
\mbox{\mbox{$\cP$}}.

If \mbox{\mbox{$(\cl{U}_{i})_{i\in I}$}} is a closed cover of \mbox{\mbox{$X$}} by trivialising
sets and \mbox{\mbox{$(\sigma_{i}\from \cl{U}_{i}\to P)_{i\in I}$}} is a collection
of continuous sections, then the collection
\mbox{\mbox{$\cl{\mathcal{U}}=(\cl{U}_{i},\sigma_{i})_{i\in I}$}}
is called a \emph{continuous closed trivialising
system} of \mbox{\mbox{$\cP$}}.

If \mbox{\mbox{$\mathcal{U}=(U_{i},\sigma_{i})_{i\in I}$}} and
\mbox{\mbox{$\cV=(V_{j},\tau _{j})_{j\in J}$}} are two continuous open
trivialising systems of \mbox{\mbox{$\cP$}}, then \mbox{\mbox{$\cV$}} is a
\emph{refinement} of \mbox{\mbox{$\mathcal{U}$}} if there exists a map
\mbox{\mbox{$J\ni j\mapsto i(j)\in I$}} such that \mbox{\mbox{$V_{j}\se U_{i(j)}$}}
and \mbox{\mbox{$\tau_{j}=\left.\sigma_{i(j)}\right|_{V_{j}}$}}, i.e.,
\mbox{\mbox{$(V_{j})_{j\in J}$}} is a refinement of \mbox{\mbox{$(U_{i})_{i\in I}$}}
and the sections \mbox{\mbox{$\tau_{j}$}} are obtained from the section
\mbox{\mbox{$\sigma_{i}$}} by restrictions.

If \mbox{\mbox{$\mathcal{U}=(U_{i},\sigma_{i})_{i\in I}$}} is a continuous open
trivialising system and \mbox{\mbox{$\cl{\cV}=(\cl{V}_{j},\tau _{j})_{j\in J}$}} is
a continuous closed trivialising system, then \mbox{\mbox{$\cl{\cV}$}} is a
\emph{refinement} of \mbox{\mbox{$\mathcal{U}$}} if there exists a map \mbox{\mbox{$J\ni
j\mapsto i(j)\in I$}} such that \mbox{\mbox{$\cl{V}_{j}\se U_{i(J)}$}} and
\mbox{\mbox{$\tau_{j}=\left.\sigma_{i(j)}\right|_{\cl{V}_{j}}$}} and vice versa.

Furthermore, if \mbox{\mbox{$\cP$}} is a smooth principal \mbox{\mbox{$K$}}-bundle
over \mbox{\mbox{$M$}}, then a \emph{smooth open trivialising system}
\mbox{\mbox{$\mathcal{U}$}} of \mbox{\mbox{$\cP$}} consists of an open cover
\mbox{\mbox{$(U_{i})_{i\in I}$}} and smooth sections \mbox{\mbox{$\sigma_{i}\from
U_{i}\to P$}}. If each \mbox{\mbox{$\cl{U}_{i}$}} is also a manifold with
corners and the section \mbox{\mbox{$\sigma_{i}$}} can be extended to smooth
sections \mbox{\mbox{$\sigma_{i}\from \cl{U}_{i}\to P$}}, then
\mbox{\mbox{$\cl{\mathcal{U}}=(\cl{U}_{i},\sigma_{i})_{i\in I}$}} is called a
\emph{smooth closed trivialising system} of \mbox{\mbox{$\cP$}}. In this
case, \mbox{\mbox{$\mathcal{U}$}} is called the trivialising system
\emph{underlying} \mbox{\mbox{$\cl{\mathcal{U}}$}}.

For each kind of trivialising system, the sections define continuous
maps \mbox{\mbox{$k_{ij}\from U_{i}\cap U_{j}\to K$}} (respectively
\mbox{$k_{ij}\from \ol{U}_{i}\cap \ol{U}_{j}\to K$} in the case of a closed
trivialising system) by
\begin{align}\label{eqn:transitionFunctions}
k_{ij}(x)=k_{\sigma_{i}}(\sigma_{j}(x))
\quad \text{ or equivalently }\quad \sigma_{i} (x)\cdot k_{ij}(x)=
\sigma_{j}(x),
\end{align}
called \emph{transition functions}.
\end{definition}

\begin{lemma}
\label{lem:forcingTransitionFunctionsIntoOpenCovers} Let \mbox{\mbox{$X$}} be a
compact space, \mbox{\mbox{$K$}} be topological group and \mbox{\mbox{$(O_{\ell})_{\ell\in L}$}}
be an open cover of \mbox{\mbox{$K$}}. If \mbox{\mbox{$\cP$}} is a continuous principal \mbox{\mbox{$K$}}-bundle
over \mbox{\mbox{$X$}}, then for each continuous open trivialising system
\mbox{\mbox{$\mathcal{U}=(U_{i},\sigma_{i})_{\In}$}} there exists a refinement
\mbox{\mbox{$\cV=(V_{s},\tau_{s})_{s=1,\dots ,r}$}} such that for each transition
function \mbox{\mbox{$k_{st}\from V_{s}\cap V_{t}\to K$}} of \mbox{\mbox{$\cV$}} we have
\mbox{\mbox{$k_{st}(V_{s}\cap V_{t})\se O_{\ell}$}} for some \mbox{\mbox{$\ell\in L$}}.
\end{lemma}

\section*{Acknowledgements}

The author of the present paper is most grateful to Professor
Karl-Hermann Neeb for the exemplary supervision of the PHD-project of
the author, out of which grew the present paper. The work on this
paper was funded by a doctoral scholarship from the Technische
Universit\"at Darmstadt.

\def\polhk#1{\setbox0=\hbox{#1}{\ooalign{\hidewidth
  \lower1.5ex\hbox{`}\hidewidth\crcr\unhbox0}}}

\vskip\baselineskip
\vskip\baselineskip
\vskip\baselineskip
\large
\noindent 
Christoph Wockel\\
Mathemtisches Institut\\
Georg-August-Universit\"at G\"ottingen\\
Bunsenstra\ss{}e 3-5\\
D-37073 G\"ottingen\\
Germany\\[\baselineskip]
\normalsize
\texttt{christoph@wockel.eu}
\end{document}